\documentclass[twoside,11pt]{article}

\usepackage[preprint]{jmlr2e}

\usepackage{enumitem}
\usepackage{graphicx,color}
\usepackage{hyperref}
\usepackage{url}
\usepackage{booktabs}
\usepackage{amsfonts}
\usepackage{nicefrac}
\usepackage{microtype}

\usepackage{amsmath, amsfonts, amssymb}
\usepackage{mathtools}
\usepackage{graphicx}
\usepackage{bm}

\newcommand{\R}{\mathbb R}
\newcommand{\X}{\mathcal X}
\newcommand{\E}{\mathbb E}

\renewcommand{\Pr}{P}
\renewcommand{\P}{\Pr}
\newcommand{\N}{\mathbb N}

\newcommand{\wh}{Y}
\newcommand{\vh}{Z}

\newcommand{\proj}{\text{proj}}
\newcommand{\argmin}{\text{argmin}}
\newcommand{\opt}{^{\star}}

\newcommand{\G}{\mathcal G}
\newcommand{\F}{\mathcal F}

\newcommand{\whp}{\text{w.p.}\geq 1-\delta}
\newcommand{\sumt}{\sum_{t=0}^{T-1}}

\newcommand{\diag}{\text{diag}}
\newcommand{\ls}{\mathcal{L}}

\newcommand{\iter}{_{\text{iter}}}
\newcommand{\emp}{_{\text{emp}}}

\usepackage{bbm}
\newcommand{\ind}{\mathbbm{I}}

\newcommand{\gtrun}[1]{\mathbb{I}_{\{ #1\} }}

\newcommand{\txi}{\widetilde{\xi}}
\newcommand{\tX}{\widetilde{X}}

\newcommand{\ev}[1]{\mathbb{E}\left[#1\right]}
\newcommand{\pr}[1]{P\left(#1\right)}

\newcommand{\vertiii}[1]{{\left\vert\kern-0.25ex\left\vert\kern-0.25ex\left\vert #1 
    \right\vert\kern-0.25ex\right\vert\kern-0.25ex\right\vert}}

\newcommand{\cmt}[1]{{\color{black}#1}}

\usepackage{lastpage}
\jmlrheading{25}{2024}{1-\pageref{LastPage}}{4/23; Revised
7/24}{8/24}{23-0466}{Liam Madden, Emiliano Dall'Anese, and Stephen Becker}

\ShortHeadings{High-probability Convergence Bounds for Non-convex SGD}{Madden, Dall'Anese, and Becker}

\begin{document}

\title{High Probability Convergence Bounds for Non-convex Stochastic Gradient Descent with Sub-Weibull Noise}

\author{\name Liam Madden \email liam@ece.ubc.ca \\
       \addr Department of Electrical and Computer Engineering\\
       University of British Columbia\\
       \AND
       \name Emiliano Dall'Anese \email edallane@bu.edu \\
       \addr Department of Electrical and Computer Engineering\\
    Boston University 
       \AND
       \name Stephen Becker \email stephen.becker@colorado.edu \\
       \addr Department of Applied Mathematics\\
       University of Colorado Boulder}

\editor{Simon Lacoste-Julien}

\maketitle

\begin{abstract}
Stochastic gradient descent is one of the most common iterative algorithms used in machine learning and its convergence analysis is a rich area of research. Understanding its convergence properties can help inform what modifications of it to use in different settings. However, 
most theoretical results either assume convexity or only provide convergence results in mean. This paper, on the other hand, proves convergence bounds in high probability without assuming convexity. Assuming strong smoothness, we prove high probability convergence bounds in two settings: (1) assuming the Polyak-Łojasiewicz inequality and norm sub-Gaussian gradient noise and (2) assuming norm sub-Weibull gradient noise. In the second setting, as an intermediate step to proving convergence, we prove a sub-Weibull martingale difference sequence self-normalized concentration inequality of independent interest. It extends Freedman-type concentration beyond the sub-exponential threshold to heavier-tailed martingale difference sequences. We also provide a post-processing method that
picks a single iterate with a 
provable 
convergence guarantee as opposed to the usual bound for the unknown best iterate.  Our convergence result for sub-Weibull noise extends the regime where stochastic gradient descent has equal or better convergence guarantees than stochastic gradient descent with modifications such as clipping, momentum, and normalization.
\end{abstract}


\section{Introduction}

Stochastic gradient descent (SGD) and its variants are some of the most commonly used algorithms in machine learning. In particular, they are used for training neural network and transformer models, models that have achieved considerable success on image classification and language processing tasks in recent years. The training in this case is non-convex and smooth for many activation/attention functions, such as sigmoid, GELU, and softmax. Even for ReLU, which is not differentiable, the training is smooth over most of the parameter space and avoidance of the non-smooth part is dealt with separately. Thus, the smooth non-convex convergence analysis of SGD has far-reaching influence on the field of machine learning.

There is a large literature on the almost sure and mean convergence of SGD assuming strong smoothness. \cite{Bertsekas} prove that SGD almost surely converges to a first-order stationary point assuming strong smoothness and the relaxed growth noise condition; more recently,
\cite{patel2020stopping} showed the same result with a different proof technique.
\cite{ghadimi2013stochastic} prove that the mean of the squared gradient norm of SGD converges to zero at a rate of $O(1/\sqrt{T})$ assuming strong smoothness and the bounded variance noise condition, where $T$ is the number of SGD iterations. \cite{Khaled} prove the same but with the expected smoothness noise condition. \citet{sebbouh2021almost} strengthen this to an almost sure convergence rate. Assuming strong smoothness and convexity, the tight mean convergence rate of the squared gradient norm of SGD is $O(1/\sqrt{T})$ \citep[Thm. 5.3.1]{Nemirovsky83}. In fact, this is the optimal rate for all stochastic first-order methods assuming only strong smoothness \citep{arjevani2019lower}. Thus, the $O(1/\sqrt{T})$ mean convergence rate of \citet{ghadimi2013stochastic} is tight and shows that SGD is optimal in the smooth non-convex setting.


However, mean convergence is not the end of the story. A convergence guarantee is generally required with some arbitrary probability, $1-\delta$.  For a single run of SGD, using Markov's inequality gives a  $1/\delta$ scaling to the bound. If one can re-run SGD many times (say, $r$ times), then by choosing the best run, $\delta$ can be relatively large since the overall success is at least by $1-\delta^r$. On the other hand, if just a single run of SGD is allowed, then the $1/\delta$ factor leads to nearly vacuous results, and hence a direct high probability convergence analysis is 
%
necessary to understand the behavior of unique runs of SGD.

\cite{li2020high} prove a $O(\log(T)\log(T/\delta)/\sqrt{T})$ high probability convergence rate for \cmt{a weighted average of the squared gradient norms of} SGD assuming strong smoothness and norm sub-Gaussian noise. \cmt{However}, a series of recent papers \citep{gurbuzbalaban2020heavy,simsekli2019tail,panigrahi19non} suggest that norm sub-Gaussian noise is often not satisfied. While the central limit theorem can be used to heuristically justify norm sub-Gaussian noise for mini-batch SGD with large batch-sizes, it cannot for small batch-sizes. In particular, \cite{panigrahi19non} provide examples where the noise is Gaussian for a batch-size of 4096, is not Gaussian for a batch-size of 32, and starts out Gaussian then becomes non-Gaussian for a batch-size of 256. \cite{gurbuzbalaban2020heavy} and \cite{simsekli2019tail} suggest instead assuming that the $p$th moment of the norm of the noise is bounded, where $p\in(1,2)$. \cite{scaman2020robustness} prove a $O(1/T^{(p-1)/p}/\delta^{(2+p)/(2p)})$ high probability convergence rate for SGD assuming strong smoothness and bounded $p$th moment noise for $p\in(1,2)$. By allowing noise with possibly infinite variance, SGD converges at a slower rate in terms of $T$. Moreover, the dependence on $1/\delta$ is polynomial rather than logarithmic. \cite{cutkosky2021high} prove a $O(\log(T/\delta)/T^{(p-1)/(3p-2)})$ high probability convergence rate for clipped SGD (which uses clipping, momentum, and normalization) assuming strong smoothness and bounded $p$th moment noise for $p\in(1,2)$, thus improving the dependence on $1/\delta$ to logarithmic, but clipped SGD will still have a slower convergence rate in terms of $T$ even when the noise is norm sub-Gaussian. So, knowledge about the type of noise should actually affect whether we use SGD or clipped SGD. This motivates the question: \textbf{How heavy can the noise be for SGD to have a high probability convergence rate that depends logarithmically on} \bm{$1/\delta$}\textbf{?}


Towards this end, we consider norm sub-Weibull noise. While it is much lighter tailed than bounded $p$th moment noise for $p\in(1,2)$, it is a natural extension of norm sub-Gaussian noise and it turns out that it admits a convergence rate for SGD with logarithmic dependence on $1/\delta$. This leads us to the first contribution of our paper. \textbf{Assuming strong smoothness and norm sub-Weibull noise, we prove a} \bm{$O((\log(T)\log(1/\delta)^{2\theta}+\log(T/\delta)^{\min\{0,\theta-1\}}\log(1/\delta))/\sqrt{T})$} \textbf{~high ~probability ~convergence ~rate ~for ~\cmt{a weighted average of the squared gradient norms of} SGD, where} \bm{$\theta$} \textbf{is the sub-Weibull tail weight.} This is given in Theorem~\ref{thm:ncconvergence}. In the special case of norm sub-Gaussian noise, which is $\theta=1/2$, this becomes $O(\log(T)\log(1/\delta)/\sqrt{T})$, which slightly improves the result of \cite{li2020high}. For more general $\theta>1/2$, we make the additional assumption that the objective function is Lipschitz continuous. At its most basic, the Lipschitz continuity assumption says that the norm of the true gradient is bounded for all of the iterates. \cite{lishaoji2022high}, building off of our pre-print, relax this to only assuming that the step-size times the true gradient is bounded for all of the iterates, which is a weaker assumption since the step-size goes to zero.

However, hidden within these convergence rates is a subtle issue that requires a post-processing algorithm. \cmt{All of these high probability convergence rates are} for a weighted average of the squared gradient norms of the iterates \cmt{rather than being for the squared gradient norm of a single iterate. While} this implies a high probability convergence rate for the iterate with the smallest squared gradient norm, to actually \cmt{determine which iterate has the smallest squared gradient norm}, we would need to estimate the true gradient for each iterate! This leads us to another question: \textbf{Is there a post-processing strategy for a single run of SGD that is as efficient as 2-RSG?}

2-RSG is the algorithm of \cite{ghadimi2013stochastic} that takes multiple runs of SGD and picks the best one to get a high probability convergence rate. First, it goes from a mean convergence rate for the iterate with the smallest squared gradient norm to a mean convergence rate for a particular iterate by randomly choosing the iterate. It does this for $\Theta(\log(1/\delta))$ runs of SGD. Then, it uses $\Theta(\log(1/\delta)\sigma^2/\delta\epsilon)$ samples to estimate the $\Theta(\log(1/\delta))$ gradients and pick the run with the smallest squared estimated gradient norm, where $\sigma^2$ is the variance of the noise and $\epsilon$ is the convergence tolerance. To compete with this, \textbf{we \cmt{introduce} a post-processing algorithm that randomly chooses} \bm{$\Theta(\log(1/\delta))$} \textbf{iterates of a single run of SGD and then picks the one with the smallest squared estimated gradient norm using} \bm{$\Theta(\log(1/\delta)\sigma^2/\delta\epsilon)$} \textbf{samples and \cmt{we prove} that our high probability convergence rate applies to the squared gradient norm of this iterate.} The algorithm \cmt{and the bound on the squared gradient norm of its output are in Theorem~\ref{thm:postprocessing}. The result can be combined with} any high probability convergence rate on a weighted average of the squared gradient norm of the iterates of SGD, such as the results of \cite{li2020high}, \cite{scaman2020robustness}, and \cite{cutkosky2021high}. Moreover, the first part of it, going from a high probability bound on a weighted average of random variables to a high probability bound on the smallest of a small subset of those random variables, applies to more general stochastic sequences \cmt{as shown in Theorem~\ref{thm:validation}, which is a probability result of independent interest.}

Underlying our convergence analysis are concentration inequalities. \cite{vladimirova2019sub}, \cite{wong2020lasso}, and \cite{bakhshizadeh2020sharp} prove concentration inequalities for sub-Weibull random variables with deterministic scale parameters. However, we need a concentration inequality for a sub-Weibull martingale difference sequence (MDS) where the scale parameters are themselves random variables making up an auxiliary sequence. MDS concentration inequalities upper bound the partial sum of the main sequence and lower bound the partial sum of the auxiliary sequence simultaneously. When the lower bound depends on the partial sum of the main sequence, we call the concentration inequality self-normalized. \cite{Freedman} proves a sub-Gaussian MDS concentration inequality. \cite{Fan} proves a sub-exponential MDS concentration inequality. \cite{Harvey1} proves a sub-Gaussian MDS self-normalized concentration inequality. The proofs for these results all rely on the moment generating function (MGF), which does not exist for sub-Weibull random variables with $\theta>1$, i.e., heavier than sub-exponential. Nevertheless, \textbf{we are able to prove a sub-Weibull martingale difference sequence self-normalized concentration inequality.} This is given in Theorem~\ref{thm:freedman}. The proof uses the MGF truncation technique of \cite{bakhshizadeh2020sharp} but applied to an MDS rather than i.i.d.\ random variables. The truncation level is determined by the tail decay rate and shows up as the $\log(T/\delta)^{\theta-1}$ in the SGD convergence rate.

We also consider convergence under the Polyak-{\L}ojasiewicz (P{\L}) inequality. Some examples of problems with P{\L} objectives are logistic regression over compact sets \citep{Karimi} and matrix factorization \citep{sun2016guaranteed}. Deep linear neural networks satisfy the P{\L} inequality in large regions of the parameter space \citep[Thm.\ 4.5]{Charles}, as do two-layer neural networks with an extra identity mapping \citep{li2017convergence}. Furthermore, sufficiently wide neural networks satisfy the P{\L} inequality locally around random initialization \citep{AZ1,AZ2,du2019shallow,Du,liu2020toward}. While strong convexity implies the P{\L} inequality, a P{\L} function need not even be convex, hence the P{\L} condition is considerably more applicable than strong convexity in the context of neural networks. Moreover, as pointed out by \cite{Karimi}, the convergence proof for gradient descent assuming strong smoothness and strong convexity is actually simplified if we use the P{\L} inequality instead of strong convexity. Thus, we would like to find a simple, elegant proof for the high probability convergence of SGD assuming strong smoothness and the P{\L} inequality.

Assuming strong smoothness and strong convexity, the tight mean convergence rate of SGD is $O(1/T)$ \citep{Nemirovski}. The same mean convergence rate can be shown assuming strong smoothness and the P{\L} inequality \citep{Karimi,Orvieto}. \cmt{But neither of these papers address high probability convergence. To fill this gap,} \textbf{we prove a} \bm{$O(\log(1/\delta)/T)$} \textbf{high probability convergence rate for the objective value of SGD assuming strong smoothness, the P}\bm{{\L}} \textbf{inequality, and norm sub-Gaussian noise.} This is given in Theorem~\ref{thm:plconvergence}. The proof relies on a novel probability result, Theorem~\ref{thm:harvey}, that essentially shows that adding two kinds of noise to a contracting sequence---sub-Gaussian noise with variance depending on the sequence itself and sub-exponential noise---results in a sub-exponential sequence. \cmt{Unfortunately, the proof does not  extend to the case of norm sub-Weibull noise. It is unclear if this is just an artifact of the proof or if the faster convergence of SGD under the P{\L} inequality really cannot be maintained under norm sub-Weibull noise.}

\cmt{
\textbf{Applicability of assumptions.} Given the motivation from neural networks, it is necessary to say a few words regarding the applicability of the Lipschitz continuity and strong smoothness assumptions. First, we show in Lemmas~\ref{lma:assump} and~\ref{lma:oymak} that these assumptions are satisfied by simple neural network models that are, nevertheless, powerful enough to interpolate generic data sets using only twice the minimal number of parameters required for \textit{any} model to do so. Second, they are satisfied if the iterates remain in a bounded set, which is precisely what happens in the neural tangent kernel setting. In this setting, the local P{\L} constant of the least squares loss is sufficiently high (due to overparameterization) for the iterates of SGD to converge to a point close to initialization~\citep{oymak2020toward}. Third, while it is true that Lipschitz continuity and strong smoothness may not be satisfied in practice, our paper is guided by the principle of ``how far can we relax assumptions while still being able to prove high probability convergence rates?'' An alternative principle is ``what can we prove given assumptions that are completely justified in practice?'' This is the principle guiding \cite{patel2022global}. They provide a neural network counterexample that does not even satisfy $(L_0,L_1)$-smoothness~\citep{zhang2020Why}, an alternative assumption to strong smoothness. Then, they prove---assuming only (1) a lower bound on the objective function, (2) local $\alpha$-H{\"o}lder continuity of the gradient, (3) equality of the true and mean gradients, and (4) that the $p$th moment of the norm of the stochastic gradient, for some $p\in(1,2]$, is bounded by an upper semi-continuous function---that, with probability 1, either (1) the norm of the iterates of SGD go to infinity, or (2) SGD almost surely converges to a first order stationary point. These assumptions are much weaker than ours, but so is the conclusion. We leave it as a future research direction to further relax our assumptions.
}

\textbf{Organization.} Section~\ref{sec:prelims} establishes the optimization framework and reviews optimization, probability, and sub-Weibull terminology and results. Section~\ref{sec:pl} proves the P{\L} convergence rate.
Section~\ref{sec:concentration} proves the MDS concentration inequality.
Section~\ref{sec:nc1} proves the non-convex convergence rate. Section~\ref{sec:postprocessing} establishes the post-processing algorithm. Section~\ref{sec:numerics} trains a two layer neural network model on synthetic data with Weibull noise injected into the gradient, showing the dependence of the tail of the convergence error on the tail of the noise.

\textbf{Notation.} We use $x$, $y$, $z$ for vectors and $X$, $Y$, $Z$, $\xi$, $\psi$ for random variables. We use $e_t$ for the error vector and distinguish it from Euler's number with the subscript (the subscript will denote the iteration index). We use $\sigma^2$ for the variance and so spell out sigma for sigma algebra. We use $O$ and $\Theta$ for big-$O$ notation. When comparing two sequences of real numbers, $(a_t)$ and $(b_t)$: $a_t=o(b_t)$ if $\lim a_t/b_t=0$, $a_t=O(b_t)$ if $\limsup |a_t|/b_t<\infty$, and $a_t=\Theta(b_t)$ if $a_t=O(b_t)$ and $b_t=O(a_t)$. We use $[n]$ to denote the set $\{1,\ldots,n\}$ and $\Gamma$ to denote the gamma function.

\section{Preliminaries}
\label{sec:prelims}

We are interested in the unconstrained optimization problem
\begin{align}
\label{eq:stochprob}
    \min_{x\in\R^d} f(x),
\end{align}
where $f: \R^d \rightarrow \R$ is a differentiable function, and the SGD iteration
\begin{align*}
    x_{t+1} = x_t-\eta_t g_t~\forall t\in\N\cup\{0\},
\end{align*}
where $g_t\in\R^d$ is an estimate of $\nabla f(x_t)$, $\eta_t$ is the step-size, and $x_0\in\R^d$ is the initial point. We restrict our attention to the setting where $x_0$ is deterministic, but the results easily extend to the setting where $x_0$ is a random vector. We define the noise as the vector $e_t=\nabla f(x_t)-g_t$ and assume that, conditioned on the previous iterates, it is unbiased and its norm is sub-Weibull (which, as we will see in Lemma~\ref{lma:moment}, implies that the variance of the noise is bounded).

One of the main examples of Eq.~\eqref{eq:stochprob} is the stochastic approximation problem: $f=\ev{F(\cdot,\xi)}$ where $(\Omega,\F,\P)$ is a probability space and $F:\R^d\times\Omega\rightarrow\R$ \citep{Nemirovski,ghadimi2013stochastic,Bottou2}. In this case, we independently sample $\xi_0,\xi_1,\ldots$ and set $g_t=\nabla F(x_t,\xi_t)$.

Another example is the sample average approximation problem which is the special case of the stochastic approximation problem where there is a finite set $\{\xi_1,\ldots,\xi_n\}$ such that $\xi=\xi_i$ with probability $1/n$. In this setting, we define $F_i=F(\cdot,\xi_i)$ so that
\begin{align*}
    f(x) = \frac{1}{n}\sum_{i=1}^n F_i(x).
\end{align*}

\subsection{Optimization}
We highlight a few key facts we use; see, e.g., \cite{Nesterov}, for more details.
All definitions are with respect to the Euclidean norm $\|\cdot\|$. Let $f:\R^d\to \R$ be differentiable, and assume $\argmin_{x\in\R^d}f(x)$ is non-empty and denote its minimum by $f\opt$.

If $f$ is continuously differentiable, then $f$ is $\rho$-Lipschitz if and only if $\|\nabla f(x)\|\leq \rho$ for all $x\in\R^d$.
We say $f$ is $L$-strongly smooth (or $L$-smooth for short) if its gradient is $L$-Lipschitz continuous. If $f$ is $L$-smooth, then a standard result using the Taylor expansion is
\begin{equation*}
f(y)\leq f(x)+\langle \nabla f(x),y-x\rangle+\frac{L}{2}\|x-y\|^2~\forall x,y\in\R^d.
\end{equation*}
Applying this result with $y=x-\frac{1}{L}\nabla f(x)$ and using $f(y)\ge f\opt$, we get
\begin{equation} \label{eq:LipschitzBound}
\|\nabla f(x)\|^2\leq 2L(f(x)-f\opt),
\end{equation}
or taking $x\in \argmin_{x'}\,f(x')$, we get
\begin{equation*}
f(y) - f^\star \le \frac{L}{2}\|y-x^\star\|^2.
\end{equation*}
Thus a convergence rate in terms of the iterates is stronger than one in terms of the objective, which is stronger than one in terms of the norm of the gradient. We say $f$ is $\mu$-Polyak-Łojasiewicz (or $\mu$-P{\L} for short) if $\|\nabla f(x)\|^2\geq 2\mu(f(x)-f\opt)~\forall x \in \R^d$. Combining this with Eq.~\eqref{eq:LipschitzBound} shows that $L\ge \mu$.

\cmt{
To justify the Lipschitz continuity and strong smoothness assumptions in the context of neural networks, we include the following lemmas, proved in Appendix~\ref{sec:assumplmapf}, which show that both properties are satisfied by the least squares loss applied to the hidden layer of a two layer neural network with, e.g., sigmoid functions such as tanh, arctan, and the logistic function  (applying the first lemma) or GELU activation (applying the second lemma). Note that such a model (with fixed outer layer), while a simple example of a neural network, is already able to interpolate generic data sets with only twice the necessary number of parameters for any model, including deeper ones, to do so~\citep{madden2024memory}.

\begin{lemma}
\label{lma:assump}
Let $m,n,d\in\N$ and $a\in\R$. Let $\phi:\R\to\R$ be twice differentiable and assume $|\phi(x)|,|\phi'(x)|,|\phi''(x)|\le a~\forall x\in\R$. Let $X\in\R^{d\times n}$, $v\in\R^m$, and $y\in\R^n$. Define
\begin{align*}
    f:\R^{d\times m}\to\R:W\mapsto \frac{1}{2}\|\phi(X^\top W)v-y\|^2.
\end{align*}
Then $f$ is Lipschitz continuous and strongly smooth.
\end{lemma}

\begin{lemma}
\label{lma:oymak}
Let $m,n,d\in\N$ and $a,b\in\R$. Let $\phi:\R\to\R$ be twice differentiable. Assume $|\phi'(x)|\le a$ and $|\phi''(x)|\le b$ for all $x\in\R$. Let $X\in\R^{d\times n}$, $v\in\R^m$, and $y\in\R^n$. Define
\begin{align*}
    f:\R^{d\times m}\to\R:W\mapsto \frac{1}{2}\|\phi(X^\top W)v-y\|^2.
\end{align*}
Let $\alpha\ge 1$ and define the sublevel set $\mathcal{W}=\{W\in\R^{d\times m}\mid f(W)\le \alpha\}$. Then, on $\mathcal{W}$, $f$ is 
\begin{align*}
    &a\|v\|_{\infty}\|X\|_2\sqrt{2\alpha m}\text{-Lipschitz continuous and}\\
    &\left(a^2\|v\|_{\infty}^2\|X\|_2^2m+b\|v\|_{\infty}\|X\|_2\|X\|_{1,2}\sqrt{2\alpha}\right)\text{-strongly smooth}.
\end{align*}
\end{lemma}
}

\subsection{Probability}

See Section 20 of \cite{Billingsley} for details on the convergence of random variables. A sequence of random variables $(X_t)$ converges to a random variable $X$:
\begin{itemize}[label={$\circ$}]
    \item in probability if, for all $\epsilon>0$, $\lim_t \P(|X_t-X|>\epsilon)=0$, denoted $X_t\overset{p}{\to}X$;
    \item in mean if $\lim_t\E|X_t-X| =0$, denoted $X_t\overset{L^1}{\to}X$;
    \item almost surely if $\P(\lim_t X_t=X)=1$, denoted $X_t\overset{a.s.}{\to}X$.
\end{itemize}
When the rate of convergence is of interest, we say a sequence of random variables $(X_t)$ converges to $X$:
\begin{itemize}[label={$\circ$}]
    \item with mean convergence rate $r(t)$ if, for all $t$, $\E|X_t-X|\leq r(t)$;
    \item with high probability convergence rate $\tilde{r}(t,\delta)$ if, for all $t$ and $\delta$, $\P\left(|X_t-X|\leq \tilde{r}(t,\delta)\,\right)\geq 1-\delta$.
\end{itemize}
All five kinds of convergence are interrelated:
\begin{itemize}[label={$\circ$}]
    \item Convergence in mean and convergence almost surely both imply convergence in probability.
    \item Mean convergence with rate $r(t)$ such that $r(t)\to 0$ as $t\to\infty$ implies convergence in mean.
    \item By Markov's inequality, mean convergence with rate $r(t)$ implies high probability convergence with rate $\tilde{r}(t,\delta)=r(t)/\delta$.
    \item By the Borel-Cantelli lemma \citep[Thm. 4.3]{Billingsley}, high probability convergence with rate $\tilde{r}(t,\delta)=r(t)p(\delta)$ implies almost sure convergence if, for all $a>0$, $\sum_{t=0}^{\infty}\min\{1,p^{-1}(a/r(t))\}<\infty$. If $p(\delta)=1/\delta^c$ for some $c>0$, then $r(t)=o(1/t^c)$ is required. If $p(\delta)=\log(1/\delta)$, then only $r(t)=O(1/t^c)$ for some $c>0$ is required.
\end{itemize}

See Section 35 of \cite{Billingsley} for details on martingale difference sequences. Let $(\Omega,\F,\P)$ be a probability space. A sequence, $(\F_i)$, of nested sigma algebras in $\F$ (i.e., $\F_i\subset \F_{i+1}\subset \F$) is called a filtration, in which case $(\Omega,\F,(\F_i),\P)$ is called a filtered probability space. A sequence of random variables $(\xi_i)$ is said to be adapted to $(\F_i)$ if each $\xi_i$ is $\F_i$-measurable. Furthermore, if $\ev{\xi_i\mid\F_{i-1}}=\xi_{i-1}~\forall i$, then $(\xi_i)$ is called a martingale. On the other hand, if $\ev{\xi_i\mid\F_{i-1}}=0~\forall i$, then $(\xi_i)$ is called a martingale difference sequence.

For the noise sequence $(e_t)$, we define a corresponding filtration $(\F_t)$ by letting $\F_t$ be the sigma algebra generated by $e_0,\ldots,e_t$ for all $t\ge 0$ and setting $\F_{-1}=\{\emptyset,\Omega\}$. Note that $x_t$ is $\F_{t-1}$-measurable while $e_t$ is $\F_t$-measurable.

\subsection{Sub-Weibull Random Variables}

We consider sub-Gaussian, sub-exponential, and sub-Weibull random variables.

\begin{definition}
A random variable $X$ is $K$-sub-Gaussian if $\ev{\exp\left(X^2/K^2\right)}\leq 2$.
\end{definition}
\begin{definition}
A random variable $X$ is $K$-sub-exponential if $\ev{\exp\left(|X|/K\right)}\leq 2$.
\end{definition}
\begin{definition}
A random variable $X$ is $K$-sub-Weibull($\theta$) if $\ev{\exp\left((|X|/K)^{1/\theta}\right)}\leq 2$.
\end{definition}

See Proposition 2.5.2 of \cite{Vershynin} for equivalent definitions of sub-Gaussian, Proposition 2.7.1 of \cite{Vershynin} for equivalent definitions of sub-exponential, and Theorem 2.1 of \cite{vladimirova2019sub} for equivalent definitions of sub-Weibull. Note that sub-Gaussian and sub-exponential are special cases of sub-Weibull using $\theta=\frac12$ and $\theta=1$, respectively. The tail parameter $\theta$ measures the heaviness of the tail---higher values correspond to heavier tails---and the scale parameter $K$ gives us the following bound on the second moment.

\begin{lemma}
\label{lma:moment}
If $X$ is $K$-sub-Weibull($\theta$) then $\ev{|X|^p}\leq 2\Gamma(\theta p+1)K^p~\forall p>0$. In particular, $\ev{X^2}\leq 2\Gamma(2\theta+1)K^2$.
\end{lemma}
\begin{proof}
First, for all $t\ge 0$,
\begin{align*}
    \pr{|X|\ge t}&=\pr{\exp\left(\left(|X|/K\right)^{1/\theta}\right)\ge \exp\left(\left(t/K\right)^{1/\theta}\right)}\\
    &\leq 2\exp\left(-\left(t/K\right)^{1/\theta}\right).
\end{align*}
Second,
\begin{align*}
    \ev{|X|^p} &= \int_0^{\infty}\pr{|X|^p\ge x}dx\\
    &\leq 2\int_0^{\infty}\exp\left(-\left(x^{1/p}/K\right)^{1/\theta}\right)\\
    &= 2\theta p K^p \int_0^{\infty}\exp(-u)u^{\theta p-1}du\\
    &= 2\theta p  \Gamma(\theta p)K^p\\
    &= 2 \Gamma(\theta p +1)K^p.
\end{align*}
\end{proof}

The proof demonstrates the general techniques for going between probability and expectation in this type of analysis. To go from probability to expectation, the CDF formula is used:
\begin{align*}
    \ev{|Y|} = \int_0^{\infty}\pr{|Y|>t}dt.
\end{align*}
To go from expectation to probability, Markov's inequality is used:
\begin{align*}
    \pr{|Y|>t}\leq \frac{1}{t}\ev{|Y|}~\forall t>0.
\end{align*}
In both cases, the trick is choosing what $Y$ should be, e.g. $Y=|X|^p$ or $Y=\exp((\lambda |X|)^{1/\theta})$.

\cmt{
\begin{remark}
\label{rmk:centering}
Note that our definitions of sub-Gaussian, sub-exponential, and sub-Weibull do not require the random variable to be centered. Thus, we do not center the constant random variable $X\equiv x$ to think of it as having scale parameter zero. Instead, for each $\theta>0$, $X\equiv x$ is $|x|/\log(2)^{\theta}$-sub-Weibull($\theta$).
\end{remark}

We can apply Remark~\ref{rmk:centering} in the following way to show how the sub-Weibull scale parameter might scale with the tail parameter for the gradient noise in the sample average approximation setting where $f(x)=\frac{1}{n}\sum_{i=1}^n F_i(x)$. Assume that $\|\nabla f(x)-\nabla F_i(x)\|\le \rho$ for all $x\in\X\subset \R^d$ and all $i\in[n]$. Since $[n]$ is a finite set, this follows for some $\rho>0$ if $\X$ is compact. Then we get that $\|\nabla f(x)-\nabla F_{\xi}(x)\|$ is $\rho/\log(2)^{\theta}$-sub-Weibull($\theta$) for all $x\in\X$ and all $\theta>0$. Thus, we can decrease the scale parameter by increasing the tail parameter. However, the convergence rate in Theorem~\ref{thm:ncconvergence} has a coefficient of $\theta^{2\theta}K^2$ when $\theta>1/2$, and $\theta^{2\theta}\rho^2/\log(2)^{2\theta}$ increases as $\theta$ increases. So, we cannot hope to erase the noise in the convergence rate by increasing the scale parameter. On the other hand, it may be the case that there is some $\theta$ slightly larger than $1/2$ which minimizes the convergence bound when all of the constants are considered.
}

\section{P{\L} Convergence}
\label{sec:pl}

In this section, we prove the following convergence bound, proving Theorem~\ref{thm:harvey}, a probability result, along the way. The convergence bound matches the optimal convergence rate in mean of $O(1/T)$ and has $\log(1/\delta)$ dependence on $\delta$.

\begin{theorem}
\label{thm:plconvergence}
Assume $f$ is $L$-smooth and $\mu$-P{\L} and that, conditioned on the previous iterates, $e_t$ is centered and $\|e_t\|$ is $\sigma$-sub-Gaussian. Then, SGD with step-size 
\begin{align*}
    \eta_t = \frac{2t+4\kappa-1}{\mu(t+2\kappa)^2},
\end{align*}
where $\kappa=L/\mu$, constructs a sequence $(x_t)$ such that, $\whp$ for all $\delta\in(0,1)$,
\begin{align*}
    f(x_T)-f\opt
    &= O\left(\frac{L\sigma^2\log(\cmt{e}/\delta)}{\mu^2 T}\right).
\end{align*}
\end{theorem}
\begin{proof}
Assume $f$ is $L$-smooth, hence 
\begin{align}
\label{eq:smoothness}
    f(y) &\leq f(x)+\langle \nabla f(x),y-x\rangle +\frac{L}{2}\|x-y\|^2
\end{align}
for all $x, y \in \R^d$. Set $y = x_{t+1}=x_t-\eta_t g_t=x_t-\eta_t(\nabla f(x_t)-e_t)$ and $x = x_{t}$ and subtract $f\opt$. Additionally assume $f$ is $\mu$-P{\L} and $\eta_t\leq 1/L$. Then
\begin{align}
    f(x_{t+1})&-f\opt\notag\\
    &\leq f(x_t)-f\opt -\eta_t\left(1-\frac{L\eta_t}{2}\right)\|\nabla f(x_t)\|^2+\eta_t(1-L\eta_t)\langle \nabla f(x_t),e_t\rangle+\frac{L\eta_t^2}{2}\|e_t\|^2\notag\\
    &\le f(x_t)-f\opt -\frac{\eta_t}{2}\|\nabla f(x_t)\|^2+\eta_t(1-L\eta_t)\langle \nabla f(x_t),e_t\rangle+\frac{L\eta_t^2}{2}\|e_t\|^2\notag\\
    &\leq (1-\mu \eta_t) (f(x_t)-f\opt)+ \eta_t(1-L\eta_t)\langle \nabla f(x_t),e_t\rangle+\frac{L\eta_t^2}{2}\|e_t\|^2\label{eq:smoothpldecrease}
\end{align}
where the second inequality used $\eta_t \leq 1/L$ and the third used the P{\L} inequality.

\cmt{Note that}
\begin{align*}
    \eta_t = \frac{2t+4\kappa-1}{\mu(t+2\kappa)^2}=\Theta\left(\frac{1}{\mu t}\right)
\end{align*}
\cmt{and,} since
$\eta_t$ is decreasing in $t$ for $t\ge 0$,
\begin{align*}
    \eta_t \leq \frac{4\kappa-1}{\mu(2\kappa)^2}\leq \frac{4\kappa}{\mu 4\kappa^2}=\frac{1}{L}, 
\end{align*}
\cmt{so} we
can apply 
Eq.~\eqref{eq:smoothpldecrease}.

Define
\begin{align*}
    &X_t = (t+2\kappa-1)^2(f(x_t)-f\opt)\\
    &\wh_t = \eta_t(1-L\eta_t)(t+2\kappa)^2\langle\nabla f(x_t),e_t\rangle\\
    &\vh_t=\frac{L\eta_t^2(t+2\kappa)^2}{2}\|e_t\|^2.
\end{align*}
Then multiplying both sides of Eq.~\eqref{eq:smoothpldecrease} by $(t+2\kappa)^2$ 
and noting
$1-\mu\eta_t = \frac{(t+2\kappa-1)^2}{(t+2\kappa)^2}$
we get the recursion 
\begin{align*}
    X_{t+1} \leq X_t+\wh_t+\vh_t.
\end{align*}

Recall that we defined $\F_t$ as the sigma algebra generated by $e_0,\ldots,e_t$ for all $t\ge 0$ and $\F_{-1}=\{\emptyset,\Omega\}$. So, $X_t$ is $\F_{t-1}$-measurable, $Y_t$ is $\F_t$-measurable, and $Z_t$ is $\F_t$-measurable.

For all $t\ge 0$, assume, conditioned on $\F_{t-1}$, that $e_t$ is centered and $\|e_t\|$ is $\sigma$-sub-Gaussian. Then, $Y_t$ is centered conditioned on $\F_{t-1}$. Note that we can bound $Y_t$ using Cauchy-Schwarz and Eq.~\eqref{eq:LipschitzBound} in the following way:
\begin{align*}
    \langle \nabla f(x_t),e_t\rangle \leq \|\nabla f(x_t)\|\cdot\|e_t\|\leq \sqrt{2L(f(x_t)-f\opt)}\|e_t\|.
\end{align*}
But, if we use this to get a new recursion, then we get a sub-optimal convergence rate. Instead, we need to keep the same recursion but use the bound on $Y_t$ in its MGF:
\begin{align*}
    \ev{\exp\left(\frac{\wh_t^2}{\frac{18L}{\mu^2}X_t\sigma^2}\right)\bigg|~\F_{t-1}} &\cmt{\le} \ev{\exp\left(\frac{\eta_t^2(1-L\eta_t)^2(t+2\kappa)^4\|\nabla f(x_t)\|^2\|e_t\|^2}{\frac{18L}{\mu^2}X_t\sigma^2}\right)\bigg|~\F_{t-1}}\\
    &\leq \ev{\exp\left(\frac{\eta_t^2(1-L\eta_t)^2\frac{(t+2\kappa)^4}{(t+2\kappa-1)^2}2LX_t\|e_t\|^2}{\frac{18L}{\mu^2}X_t\sigma^2}\right)\bigg|~\F_{t-1}}
    \text{ via }\eqref{eq:LipschitzBound}
    \\
    &\leq \ev{\exp\left(\frac{\frac{18L}{\mu^2}X_t\|e_t\|^2}{\frac{18L}{\mu^2}X_t\sigma^2}\right)\bigg|~\F_{t-1}}\\
    &\leq 2.
\end{align*}
Thus, $\wh_t$ is $\frac{18L}{\mu^2}X_t\sigma^2$-sub-Gaussian conditioned on $\F_{t-1}$. Similarly,
\begin{align*}
    \ev{\exp\left(\frac{\vh_t}{\frac{2L}{\mu^2}\sigma^2}\right)\bigg|~\F_{t-1}} &\leq \ev{\exp\left(\frac{\frac{2L}{\mu^2}\|e_t\|^2}{\frac{2L}{\mu^2}\sigma^2}\right)\bigg|~\F_{t-1}}\\
    &\leq 2.
\end{align*}
and so $\vh_t$ is $\frac{2L}{\mu^2}\sigma^2$-sub-exponential conditioned on $\F_{t-1}$.

So, we have a recursion for a stochastic sequence with two types of additive noise---sub-Gaussian noise, with variance depending on the sequence itself, and sub-exponential noise. We prove that this makes the main sequence sub-exponential in the following theorem.

\begin{theorem}
\label{thm:harvey}
Let $(\Omega,\F,(\F_i),P)$ be a filtered probability space. Let $\cmt{(X_{i+1})}$, $(\wh_i)$, and $(\vh_i)$ be adapted to $(\F_i)$, and $X_0$ be deterministic. \cmt{Let $(\alpha_i)$ and $(\beta_i)$ be sequences of non-negative reals and let $(\gamma_i)$ be a sequence of positive reals.} Assume $X_i$ and $\vh_i$ are non-negative almost surely. Assume $\E[\exp(\lambda \wh_i) \mid \F_{i-1}]\leq \exp(\frac{\lambda^2}{2}\beta_i^2 X_i)$ for all $\lambda\in\R$ and $\E[\exp(\lambda \vh_i) \mid \F_{i-1}]\leq \exp(\lambda \gamma_i)$ for all $\lambda\in [0,\gamma_i^{-1}]$. Assume
\begin{align*}
    X_{i+1}\leq \alpha_i X_i+\wh_i+\vh_i.
\end{align*}
Then, for any \cmt{sequence of positive reals} $(K_i)$ such that $K_0\geq X_0$ and, for all $i\geq 0$, $K_{i+1}^2\geq (\alpha_i K_i+2\gamma_i)K_{i+1}+\beta_i^2 K_i$, and \cmt{for }any $n\ge0$, we have, $\whp$ for all $\delta \in (0,1)$,
\begin{align*}
    X_n \leq K_n\log(e/\delta).
\end{align*}
\end{theorem}
\begin{proof}
We want to find requirements on a sequence $(K_i)$ such that $\ev{\exp(\lambda X_i)}\leq\exp(\lambda K_i)~\forall \lambda\in[0,K_i^{-1}]$. Then, by Markov's inequality and taking $\lambda = K_{n}^{-1}$,\\ $\pr{X_n\ge \log(e/\delta)K_n}=\pr{\exp(X_n/K_n)\ge e/\delta}\leq \delta$. Our proof is inductive. For the base case, we need
\begin{align*}
    (1)~K_0\ge X_0.
\end{align*}
For the induction step, assume $\ev{\exp(\tilde{\lambda} X_i)}\leq\exp(\tilde\lambda K_i)~\forall\tilde\lambda\in[0,K_i^{-1}]$. Let $\lambda\in[0,K_{i+1}^{-1}]$. Then
\begin{align*}
    \ev{\exp(\lambda X_{i+1})} &\leq \ev{\exp(\lambda\alpha_iX_i+\lambda\wh_i+\lambda \vh_i)}\quad\text{by non-negativity}\\
    &\overset{\text{TE}}{=}\ev{\exp(\lambda\alpha_iX_i)\ev{\exp(\lambda\wh_i)\exp(\lambda \vh_i)\mid\F_{i-1}}}\\
    &\overset{\text{CS}}{\leq}\ev{\exp(\lambda\alpha_iX_i)\ev{\exp(2\lambda\wh_i)\mid\F_{i-1}}^{1/2}\ev{\exp(2\lambda \vh_i)\mid\F_{i-1}}^{1/2}}
     \\
    &\overset{(2)}{\leq}\ev{\exp(\lambda \alpha_i X_i)\exp\left(\frac{(2\lambda)^2}{2}\beta_i^2 X_i\right)^{1/2}\exp(2\lambda \gamma_i)^{1/2}} \quad\text{if}\quad  2\lambda\in[0,\gamma_i^{-1} ]\\
    &=\ev{\exp(\lambda\alpha_iX_i+\lambda^2\beta_i^2X_i+\lambda\gamma_i)}\\
    &=\ev{\exp(\lambda(\alpha_i+\lambda\beta_i^2)X_i}\exp(\lambda\gamma_i)\\
    &\overset{(3)}{\leq}\exp(\lambda((\alpha_i+\lambda\beta_i^2)K_i+\gamma_i))
    \quad\text{via induction if}\; \tilde\lambda := \lambda(\alpha_i+\lambda\beta_i^2) \le K_i^{-1}
    \\
    &\overset{(4)}{\leq}\exp(\lambda K_{i+1})
\end{align*}
where $TE$ denotes the law of total expectation, $CS$ denotes the Cauchy-Schwarz inequality, and $(2)-(4)$ are the requirements
\begin{align*}
    &(2)~2\lambda\leq\gamma_i^{-1}\impliedby 2K_{i+1}^{-1}\leq \gamma_i^{-1}\iff K_{i+1}\ge 2\gamma_i\\
    &(3)~\lambda(\alpha_i+\lambda\beta_i^2)\leq K_{i}^{-1} \impliedby K_{i+1}^{-1}(\alpha_i+K_{i+1}^{-1}\beta_i^2)\leq K_{i}^{-1} \iff K_{i+1}^2\ge \alpha_iK_iK_{i+1}+\beta_i^2K_i\\
    &(4)~K_{i+1}\ge (\alpha_i+\lambda\beta_i^2)K_i+\gamma_i \impliedby K_{i+1}\ge (\alpha_i+K_{i+1}^{-1}\beta_i^2)K_i+\gamma_i\\
    &\hspace{7cm}\iff K_{i+1}^2\ge(\alpha_iK_i+\gamma_i)K_{i+1}+\beta_i^2K_i.
\end{align*}
The assumptions of the theorem imply requirements $(2)-(4)$, completing the proof.
\end{proof}

The proof is similar to the proof of Theorem 4.1 of \cite{Harvey1} but with some key differences. There, the recursion is $X_{i+1}\leq \alpha_i X_i+\wh_i \sqrt{X_i}+\gamma_i$ where $\ev{\exp(\lambda \wh_i)\mid \F_{i-1}}\leq \exp\left(\frac{\lambda^2}{2}\beta_i^2\right)$ and $\gamma_i$ is deterministic. We had to move the implicit dependence of the sub-Gaussian term inside of the MGF. We also had to allow $\gamma_i$ to be a sub-exponential random variable $\vh_i$, and so applied Cauchy-Schwarz, which contributed the 2 in the recursion for $(K_i)$.

\cmt{Note that if $\alpha_i=1,\gamma_i=\gamma,\beta_i=\beta~\forall i$, then the assumptions on $(K_i)$ are satisfied by, for example, $K_0=X_0$ and $K_{i+1}=K_i+2\gamma+\beta^2~\forall i$. Returning to the proof of Theorem~\ref{thm:plconvergence}}, for absolute constants $a$ and $b$, by Proposition 2.5.2 and Proposition 2.7.1 of \cite{Vershynin}, we can set $\beta_t^2=\frac{18aL\sigma^2}{\mu^2}$ and $\gamma_t=\frac{2bL\sigma^2}{\mu^2}$ and apply Theorem~\ref{thm:harvey} with
\begin{align*}
    K_T &\cmt{\coloneqq} X_0+\sumt (2\gamma_t+\beta_t^2)\\
    &= X_0+\frac{(18a+\cmt{4}b)L\sigma^2T}{\mu^2}.
\end{align*}
Dividing by $(T+2\kappa)^2$ completes the proof.
\end{proof}

\begin{remark}
It is natural to ask whether we can relax the sub-Gaussian assumption to a sub-Weibull assumption. In Theorem~\ref{thm:harvey}, we need a bound on the MGFs of both $\wh_i$ and $\vh_i$. But, if $\|e_t\|$ is $\sigma$-sub-Weibull($\theta$) with $\theta>1/2$, then $\|e_t\|^2$ is $\sigma^2$-sub-Weibull($2\theta$). Thus, $\vh_i$ would be sub-Weibull with tail parameter greater than 1, and so may have an infinite moment generating function for all $\lambda> 0$.
\end{remark}

A big take-away from the P{\L} analysis is its simplicity, matching the simplicity of gradient descent's convergence analysis in the same setting. It also serves as a good warm-up for the non-convex analysis that follows. Induction on the convergence recursion does not work for the non-convex analysis. Instead it requires an MDS self-normalized concentration inequality, \cmt{which we state and prove in the next section.}

\section{Concentration Inequality}
\label{sec:concentration}

In this section, we prove the following MDS concentration inequality. For $\theta=1/2$ without $\alpha$ (i.e., $\alpha=0$), Eq.~\eqref{eq:freedmansubgnoalpha}, we recover the classical Freedman's inequality \citep{Freedman}. For $\theta\in(1/2,1]$ without $\alpha$, Eq.~\eqref{eq:freedmansubenoalpha}, we recover Theorem~2.6 of \cite{Fan}. For $\theta=1/2$ with $\alpha$, Eq.~\eqref{eq:freedmansubg}, we recover Theorem~3.3 of \cite{Harvey1}, called the ``Generalized Freedman'' inequality.

\begin{theorem}
\label{thm:freedman}
Let $(\Omega,\F,(\F_i),\P)$ be a filtered probability space. Let $(\xi_i)$ and $(K_i)$ be adapted to $(\F_i)$. Let $n\in\N$. For all $i\in[n]$, assume $K_{i-1}\ge 0$ almost surely, $\ev{\xi_i\mid\F_{i-1}}=0$, and
\begin{align*}
    \ev{\exp\left((|\xi_i|/K_{i-1})^{1/\theta}\right)\mid\F_{i-1}}\leq 2
\end{align*}
where $\theta\ge 1/2$. If $\theta>1/2$, assume there exist constants $(m_i)$ such that $K_{i-1}\leq m_i$ almost surely for all $i\in[n]$.

If $\theta=1/2$, then for all $x,\beta\ge 0$, and $\alpha>0$, and $\lambda\in\left[0,\frac{1}{2\alpha}\right]$,
\begin{equation}
\label{eq:freedmansubg}
    \pr{\bigcup_{k\in[n]}\bigg\{\sum_{i=1}^k\xi_i\ge x\text{ and }\sum_{i=1}^k 2K_{i-1}^2\leq \alpha\sum_{i=1}^k \xi_i+\beta\bigg\}}\leq \exp(-\lambda x+2\lambda^2\beta),
\end{equation}
and for all $x,\beta,\lambda\ge 0$,
\begin{equation}
\label{eq:freedmansubgnoalpha}
    \pr{\bigcup_{k\in[n]}\bigg\{\sum_{i=1}^k\xi_i\ge x\text{ and }\sum_{i=1}^k 2K_{i-1}^2\leq \beta\bigg\}}\leq \exp\left(-\lambda x+\frac{\lambda^2}{2}\beta\right).
\end{equation}

If $\theta\in\left(\frac{1}{2},1\right]$, define
\begin{align*}
    &a = (4\theta)^{2\theta}e^2\\
    &b = (4\theta)^{\theta}e.
\end{align*}
For all $x,\beta\ge 0$, and $\alpha\ge b \max_{i\in[n]} m_i$, and $\lambda\in\left[0,\frac{1}{2\alpha}\right]$,
\begin{equation}
\label{eq:freedmansube}
    \pr{\bigcup_{k\in[n]}\bigg\{\sum_{i=1}^k\xi_i\ge x\text{ and }\sum_{i=1}^k aK_{i-1}^2\leq \alpha\sum_{i=1}^k \xi_i+\beta\bigg\}}\leq \exp(-\lambda x+2\lambda^2\beta),
\end{equation}
and for all $x,\beta\ge 0$, and $\lambda \in\left[0,\frac{1}{b\max_{i\in[n]}m_i}\right]$,
\begin{equation}
\label{eq:freedmansubenoalpha}
    \pr{\bigcup_{k\in[n]}\bigg\{\sum_{i=1}^k\xi_i\ge x\text{ and }\sum_{i=1}^k aK_{i-1}^2\leq \beta\bigg\}}\leq \exp\left(-\lambda x+\frac{\lambda^2}{2}\beta\right).
\end{equation}

If $\theta>1$, let $\delta\in(0,1)$. Define
\begin{align*}
    &a = (2^{2\theta+1}+2)\Gamma(2\theta+1)+\frac{2^{3\theta}\Gamma(3\theta+1)}{3\log(n/\delta)^{\theta-1}}\\
    &b = 2\log(n/\delta)^{\theta-1}.
\end{align*}
For all $x,\beta\ge 0$, and $\alpha \ge b \max_{i\in[n]}m_i$, and $\lambda\in\left[0,\frac{1}{2\alpha}\right]$,
\begin{align*}
    \pr{\bigcup_{k\in[n]}\bigg\{\sum_{i=1}^k\xi_i\ge x\text{ and }\sum_{i=1}^k aK_{i-1}^2\leq \alpha\sum_{i=1}^k \xi_i+\beta\bigg\}}&\leq \exp(-\lambda x+2 \lambda^2\beta)+2\delta,
\end{align*}
and for all $x,\beta\ge 0$, and $\lambda \in\left[0,\frac{1}{b\max_{i\in[n]}m_i}\right]$,
\begin{equation}
\label{eq:freedmansubwnoalpha}
    \pr{\bigcup_{k\in[n]}\bigg\{\sum_{i=1}^k\xi_i\ge x\text{ and }\sum_{i=1}^k aK_{i-1}^2\leq \beta\bigg\}}\leq \exp\left(-\lambda x+\frac{\lambda^2}{2}\beta\right)+2\delta.
\end{equation}
\end{theorem}
\begin{proof}
Let $(\Omega,\F,(\F_i),\P)$ be a filtered probability space. Let $(\xi_i)$ and $(K_i)$ be adapted to $(\F_i)$. Let $n\in\N$. For all $i\in[n]$, assume $0\leq K_{i-1}\leq m_i$ almost surely, $\ev{\xi_i\mid\F_{i-1}}=0$, and
\begin{align*}
    \ev{\exp\left((|\xi_i|/K_{i-1})^{1/\theta}\right)\mid\F_{i-1}}\leq 2.
\end{align*}
What we want to upper bound in this setting is
\begin{align*}
    \pr{\bigcup_{k\in[n]}\bigg\{\sum_{i=1}^k\xi_i\ge x\text{ and }\sum_{i=1}^k K_{i-1}^2\leq \alpha\sum_{i=1}^k \xi_i+\beta\bigg\}}
\end{align*}
for constants $x,\alpha>0$ and $\beta \ge 0$. To understand how to use this, set $\beta=0$ and observe
\begin{align*}
    &\pr{\sum_{i=1}^n\xi_i\ge x+\frac{1}{\alpha}\sum_{i=1}^n K_{i-1}^2}\\
    &\leq\pr{\sum_{i=1}^n\xi_i\ge x\text{ and }\sum_{i=1}^n K_{i-1}^2\leq\alpha \sum_{i=1}^n\xi_i}\\
    &\leq\pr{\bigcup_{k\in[n]}\bigg\{\sum_{i=1}^k\xi_i\ge x\text{ and }\sum_{i=1}^k K_{i-1}^2\leq\alpha \sum_{i=1}^k\xi_i\bigg\}}
\end{align*}
by monotonicity ($A\subset B\implies P(A)\leq P(B)$). The stronger bound over the union of partial sum bounds does not help us since the constants cannot depend on $k$. The union of partial sum bounds arises in the proof from a stopping time defined as the first $k$ such that the corresponding set occurs.

Proving such a bound would typically involve using an MGF bound for $\xi_i$, but the MGF is infinite in our setting. To get around this, we truncate $\xi_i$ to an appropriate level. By accounting for the probability of exceeding truncation, applying an MGF bound for truncated random variables, and applying a concentration inequality for bounded MGF martingale different sequence, we are able to prove a concentration inequality for sub-Weibull martingale difference sequences.

\textbf{Probability of exceeding truncation.} Let $\delta\in (0,1)$. We want to define the truncation level so that the probability of any of the $\xi_i$ exceeding it is smaller than $O(\delta)$. This ends up being $\txi_i = \xi_i\gtrun{\xi_i\leq cK_{i-1}}$ with $c = \log(n/\delta)^{\theta}$ as we will see.

First, for any $c>0$,
\begin{align*}
    \pr{|\xi_i|\ge c K_{i-1}} &= \pr{\exp\left((|\xi_i|/K_{i-1})^{1/\theta}\right)\ge \exp\left(c^{1/\theta}\right)}\\
    &\leq \exp\left(-c^{1/\theta}\right)\ev{\exp\left((|\xi_i|/K_{i-1})^{1/\theta}\right)}\\
    &= \exp\left(-c^{1/\theta}\right)\ev{\ev{\exp\left((|\xi_i|/K_{i-1})^{1/\theta}\right)\mid\F_{i-1}}}\\
    &\leq 2\exp\left(-c^{1/\theta}\right)
\end{align*}
where we use Markov's inequality, the law of total expectation, and the sub-Weibull assumption. In particular, $\pr{\xi_i > \log(n/\delta)^{\theta}K_{i-1}}\leq 2\delta/n$.

Then, we can bound
\begin{align*}
    \pr{\bigcup_{k\in[n]}\bigg\{\sum_{i=1}^k\xi_i\ge x\text{ and }\sum_{i=1}^k K_{i-1}^2\leq \alpha\sum_{i=1}^k \xi_i+\beta\bigg\}}
\end{align*}
by
\begin{align*}
    \pr{\bigcup_{k\in[n]}\bigg\{\sum_{i=1}^k\txi_i\ge x\text{ and }\sum_{i=1}^k K_{i-1}^2\leq \alpha\sum_{i=1}^k \txi_i+\beta\bigg\}}
\end{align*}
and
\begin{align*}
    \pr{\bigcup_{i\in[n]}\bigg\{\xi_i>\log(n/\delta)^{\theta}K_{i-1}\bigg\}} &\leq \sum_{i=1}^n\pr{\xi_i>\log(n/\delta)^{\theta}K_{i-1}}\\
    &\leq 2\delta.
\end{align*}
and so proceed using $\txi_i$ while carrying around an additional $2\delta$ probability.

\textbf{MGF bound of truncated random variable.} In order to bound the MGF of our truncated random variables, we prove the following lemma which slightly modifies Corollary 2 of \cite{bakhshizadeh2020sharp} using the law of total expectation. The main subtlety in the extension is where we use the following bound
\begin{align*}
    \pr{|X|\ge t\mid\G}&= \pr{\exp\left((|X|/K_0)^{1/\theta}\right)\ge \exp\left((t/K_0)^{1/\theta}\right)\mid \G}\\
    &\leq \exp\left(-(t/K_0)^{1/\theta}\right)\ev{\exp\left((|X|/K_0)^{1/\theta}\right)\mid \G}\\
    &\leq 2\exp\left(-(t/K_0)^{1/\theta}\right)
\end{align*}
which we denote by $*$. Otherwise, the proof exactly follows theirs.

\begin{lemma}
\label{lma:bakhshizadeh}
Let $(\Omega,\F,\P)$ be a probability space, $\G\subset\F$ be a sigma algebra, and $X$ and $K_0$ be random variables. Assume $K_0$ is $\G$-measurable. Assume, conditioned on $\G$, that $X$ is centered and $K_0$-sub-Weibull($\theta$) with $\theta>1$. Define $\tX=X\gtrun{X\leq cK_0}$. Then
\begin{align*}
    \ev{\exp(\lambda \tX)\mid \G}\leq \exp\left(\frac{\lambda^2}{2}aK_0^2\right)~\forall\lambda\in\left[0,\frac{1}{bK_0}\right]
\end{align*}
where
\begin{align*}
    &a = (2^{2\theta+1}+2)\Gamma(2\theta+1)+\frac{2^{3\theta}\Gamma(3\theta+1)}{3c^{1-1/\theta}}\\
    &b = 2c^{1-1/\theta}.
\end{align*}
\end{lemma}
\begin{proof}
For convenience, define $L_0=cK_0$. Let $\lambda\in\left[0,\frac{1}{bK_0}\right]$. That is, $\lambda\ge 0$ and $\lambda L_0^{1-1/\theta}K_0^{1/\theta}\leq \frac{1}{2}$. Since $\lambda\ge 0$, we have, by Lemma 4 of \cite{bakhshizadeh2020sharp},
\begin{align*}
    \ev{\exp(\lambda \tX)\mid\G}\leq \exp\bigg(\frac{\lambda^2}{2}\bigg(\ev{X^2\gtrun{X<0}\mid\G}+\ev{X^2\exp(\lambda X)\gtrun{0\leq X\leq L_0}\mid\G}\bigg)\bigg).
\end{align*}
Observe
\begin{align*}
    \ev{X^2\gtrun{X<0}\mid\G} &= \int_0^{\infty}\pr{X^2\gtrun{X<0}>x\mid\G}dx\\
    &= \int_0^{\infty}\pr{X^2>t^2,X<0\mid\G}2t\,dt\\
    &\leq \int_0^{\infty}\pr{|X|> t\mid\G}2t\,dt\\
    &\overset{*}{\leq} 2\int_0^{\infty}\exp\left(-(t/K_0)^{1/\theta}\right)2t\,dt\\
    &= 2\Gamma(2\theta+1)K_0^2
\end{align*}
and
\begin{align*}
    &\ev{X^2\exp(\lambda X)\gtrun{0\leq X\leq L_0}\mid\G}\\
    &= \int_0^{\infty}\pr{X^2\exp(\lambda X)\gtrun{0\leq X\leq L_0}>x\mid\G}dx\\
    &=\int_0^{\infty}\pr{X^2\exp(\lambda X)>t^2\exp(\lambda t),0\leq X\leq L_0\mid\G}(2t+\lambda t^2)\exp(\lambda t)dt\\
    &= \int_0^{\infty}\pr{|X|>t,0\leq X\leq L_0\mid\G}(2t+\lambda t^2)\exp(\lambda t)dt\\
    &\leq\int_0^{L_0}\pr{|X|>t\mid\G}(2t+\lambda t^2)\exp(\lambda t)dt\\
    &\overset{*}{\leq} 2 \int_0^{L_0} \exp\left(-\left(1-\lambda t^{1-1/\theta}K_0^{1/\theta}\right)(t/K_0)^{1/\theta}\right)(2t+\lambda t^2)dt\\
    &\leq 2 \int_0^{L_0} \exp\left(-\left(1-\lambda L_0^{1-1/\theta}K_0^{1/\theta}\right)(t/K_0)^{1/\theta}\right)(2t+\lambda t^2)dt\\
    &=2\int_0^{L_0}\exp(-u)\left(\frac{2K_0^2\theta u^{2\theta-1}}{\left(1-\lambda L_0^{1-1/\theta}K_0^{1/\theta}\right)^{2\theta}}+\frac{\lambda K_0^3 \theta u^{3\theta-1}}{\left(1-\lambda L_0^{1-1/\theta}K_0^{1/\theta}\right)^{3\theta}}\right)du\\
    &=\frac{2K_0^2\Gamma(2\theta+1)}{\left(1-\lambda L_0^{1-1/\theta}K_0^{1/\theta}\right)^{2\theta}}+\frac{2\lambda K_0^3\Gamma(3\theta+1)}{3\left(1-\lambda L_0^{1-1/\theta}K_0^{1/\theta}\right)^{3\theta}}\\
    &\leq \left(2^{2\theta+1}\Gamma(2\theta+1)+\frac{2^{3\theta}K_0\Gamma(3\theta+1)}{3L_0^{1-1/\theta}K_0^{1/\theta}}\right)K_0^2\\
    &=\left(2^{2\theta+1}\Gamma(2\theta+1)+\frac{2^{3\theta}\Gamma(3\theta+1)}{3\log(n/\delta)^{\theta-1}}\right)K_0^2.
\end{align*}
\end{proof}

Applying this to our setting, we get the MGF bound
\begin{align}
\label{eq:trunmgf}
    \ev{\exp(\lambda \txi_i)\mid \F_{t-1}}\leq \exp\left(\frac{\lambda^2}{2}aK_{i-1}^2\right)~\forall\lambda\in\left[0,\frac{1}{bK_{i-1}}\right]
\end{align}
where
\begin{align*}
    &a = (2^{2\theta+1}+2)\Gamma(2\theta+1)+\frac{2^{3\theta}\Gamma(3\theta+1)}{3\log(n/\delta)^{\theta-1}}\\
    &b = 2\log(n/\delta)^{\theta-1}.
\end{align*}

\textbf{Concentration inequality for MGF bound.} Now we just need a self-normalized concentration inequality in the following setting:

Let $(\Omega,\F,(\F_i),\P)$ be a filtered probability space. Let $(\xi_i)$ and $(K_i)$ be adapted to $(F_i)$. Let $n\in\N$. For all $i\in[n]$, assume $0\leq K_{i-1}\leq m_i$ almost surely, $\ev{\xi\mid\F_{i-1}}= 0$, and
\begin{align*}
    \ev{\exp(\lambda \xi_i)\mid\F_{i-1}} \leq \exp\left(\frac{\lambda^2}{2}aK_{i-1}^2\right)~\forall \lambda\in\left[0,\frac{1}{bK_{i-1}}\right]
\end{align*}
for constants $a,b>0$.

We can recognize this as a centered sub-exponential type MGF bound from Proposition 2.7.1 of \cite{Vershynin}. Theorem 2.6 of \cite{Fan} proves a concentration inequality in this setting, but not a self-normalized one. On the other hand, Theorem 3.3 of \cite{Harvey1} proves a self-normalized concentration inequality in the centered sub-Gaussian setting, thus extending the original result of \cite{Freedman} to a self-normalized result. We use the same trick of \cite{Harvey1} to extend Theorem 2.6 of \cite{Fan} to a self-normalized result.



We generalize the proof of Theorem 2.1 in \cite{Fan} so that it can be applied more generally, including to our setting.




\begin{lemma}
\label{lma:fan}
Let $(\Omega,\F,(\F_i),P)$ be a filtered probability space. If $(\psi_i)$ is adapted to $(\F_i)$, each $\psi_i$ is almost surely nonzero, and $(A_k)$ is a sequence of events; then
\begin{align*}
    \pr{\bigcup_{k=1}^nA_k}\le \sup_{\omega\in\Omega}\sup_{k\in[n]}\bigg\{\ind_{A_k}\prod_{i=1}^k\frac{\ev{\psi_i\mid \F_{i-1}}}{\psi_i}(\omega)\bigg\}.
\end{align*}
\end{lemma}
\begin{proof}
Define $Z_k=\prod_{i=1}^k\frac{\psi_i}{\ev{\psi_i\mid \F_{i-1}}}$. By construction, $(Z_k)$ is a martingale with initial value 1, as is $(Z_{k\land T})$ for any stopping time $T$ (where $a\land b$ denotes $\min\{a,b\}$), so $\ev{Z_{n\land T}}=1$ for any stopping time $T$. Define $T$ by $T(\omega)=\min\{k\in[n]\mid \omega\in A_k\}$.
Note that $\ind_{\cup_{k\in[n]} A_k}=\sum_{i=1}^n\ind_{\{T=k\}}$.
Thus, by the Cauchy-Schwarz inequality,
\begin{align*}
    \sum_{k=1}^n \ev{Z_k\ind_{\{T=k\}}}= \ev{Z_{n\land T}\ind_{\cup_{k\in[n]}A_k}}\le \ev{Z_{n\land T}}^{1/2}\ev{\ind_{\cup_{k\in[n]}A_k}}^{1/2}= \pr{\bigcup_{k=1}^nA_k}^{1/2}.
\end{align*}
Also,
\begin{align*}
    \pr{\bigcup_{k=1}^nA_k} &= \ev{\ind_{\cup_{k\in[n]}A_k}}=\ev{\sum_{k=1}^n\ind_{\{T=k\}}}=\sum_{k=1}^n\ev{\ind_{\{T=k\}}}.
\end{align*}
So, multiplying by $\ind_{A_k} Z_k^{-1}Z_k$ and applying Hölder's inequality,
\begin{align*}
    \pr{\bigcup_{k=1}^nA_k} &=\sum_{k=1}^n\ev{\ind_{A_k}Z_k^{-1}Z_k\ind_{\{T=k\}}}\\
    &\le \sum_{k=1}^n \sup_{\omega\in\Omega}\bigg\{\ind_{A_k}Z_k^{-1}(\omega)\bigg\}\ev{Z_k\ind_{\{T=k\}}}\\
    &\le \sup_{\omega\in\Omega}\sup_{k\in[n]}\bigg\{\ind_{A_k}Z_k^{-1}(\omega)\bigg\}\sum_{k=1}^n \ev{Z_k\ind_{\{T=k\}}}\\
    &\le \sup_{\omega\in\Omega}\sup_{k\in[n]}\bigg\{\ind_{A_k}Z_k^{-1}(\omega)\bigg\} \pr{\bigcup_{k=1}^nA_k}^{1/2}.
\end{align*}
Using $\pr{\cup_{k\in[n] A_k}}\le 1$ proves the lemma.
\end{proof}


We apply Lemma~\ref{lma:fan} to our setting to get the following self-normalized concentration inequality.

\begin{lemma}
\label{lma:subexpharvey}
Let $(\Omega,\F,(\F_i),\P)$ be a filtered probability space. Let $(\xi_i)$ and $(K_i)$ be adapted to $(F_i)$. Let $n\in\N$. For all $i\in[n]$, assume $0\leq K_{i-1}\leq m_i$ almost surely, $\ev{\xi\mid\F_{i-1}}= 0$, and
\begin{align*}
    \ev{\exp(\lambda \xi_i)\mid\F_{i-1}} \leq \exp\left(\frac{\lambda^2}{2}aK_{i-1}^2\right)~\forall \lambda\in\left[0,\frac{1}{bK_{i-1}}\right]
\end{align*}
for constants $a,b>0$. Then, for all $x,\beta\ge 0$, and $\alpha \ge b \max_{i\in[n]}m_i$, and $\lambda\in\left[0,\frac{1}{2\alpha}\right]$,
\begin{align*}
    \pr{\bigcup_{k\in[n]}\bigg\{\sum_{i=1}^k\xi_i\ge x\text{ and }\sum_{i=1}^k aK_{i-1}^2\leq \alpha\sum_{i=1}^k \xi_i+\beta\bigg\}}&\leq \exp(-\lambda x+2 \lambda^2\beta).
\end{align*}
\end{lemma}
\begin{proof}
By Claim C.2 of \cite{Harvey1}, if $0\leq\lambda\leq \frac{1}{2\alpha}$, then $\exists~ c\in[0,2]$ such that $\frac{1}{2}(\lambda+\alpha c\lambda^2)^2=c\lambda^2$. Define
\begin{align*}
    \psi_i = \exp\left((\lambda+\alpha c\lambda^2)\xi_i\right).
\end{align*}
With $0\leq\lambda\leq\frac{1}{2\alpha}$, we want $\lambda+\alpha c\lambda^2\leq\frac{1}{bK_{i-1}}$. This is ensured by $\alpha \ge b \max_{i\in[n]}m_i$.

Define
\begin{align*}
    A_k =\bigg\{\sum_{i=1}^k \xi_i\ge x\text{ and }\sum_{i=1}^k aK_{i-1}^2\leq \alpha\sum_{i=1}^k \xi_i+\beta\bigg\}.
\end{align*}
Then $\omega\in A_k$ implies
\begin{align*}
    \prod_{i=1}^k \frac{\ev{\psi_i\mid\F_{i-1}}}{\psi_i} &\leq \exp\left(-(\lambda+\alpha c\lambda^2)\sum_{i=1}^k \xi_i+\frac{(\lambda+\alpha c\lambda^2)^2}{2}\sum_{i=1}^k aK_{i-1}^2\right)\\
    &\leq \exp(-\lambda x+c\lambda^2\beta)\\
    &\leq \exp(-\lambda x+2\lambda^2\beta).
\end{align*}
\end{proof}

\textbf{Final step.} Putting everything together proves the $\theta>1$ and $\alpha>0$ case. The rest of the work for the other cases is included in Appendix~\ref{sec:freedmanpf}.
\end{proof}

\section{Non-convex Convergence}
\label{sec:nc1}

In this section, we prove the following convergence bound.

\begin{theorem}
\label{thm:ncconvergence}
Assume $f$ is $L$-smooth and that, conditioned on the previous iterates, $e_t$ is centered and $\|e_t\|$ is $K\text{-sub-Weibull}(\theta)$ with $\theta\ge 1/2$. If $\theta>1/2$, assume $f$ is $\rho$-Lipschitz. Let $\delta_1,\delta_2,\delta_3\in(0,1)$ and define $\delta=\max\{\delta_1,\delta_2,\delta_3\}$. Then, for $T$ iterations of SGD with $\eta_t=c/\sqrt{t+1}$ where $c\leq 1/L$, w.p.~$\ge 1-\delta_1-\delta_2-\delta_3$,
\begin{align*}
    &\frac{1}{\sqrt{T}}\sumt \frac{1}{\sqrt{t+1}}\|\nabla f(x_t)\|^2\\
    &\hspace{1cm}\leq \frac{4(f(x_0)-f\opt)}{c\sqrt{T}}+\frac{\gamma(\theta)\log(1/\delta_3)}{\sqrt{T}}+ \frac{4LK^2c(4e\theta\log(2/\delta_1))^{2\theta}\log(T+1)}{\sqrt{T}}\\
    &\hspace{1cm}= O\left(\frac{\log(T)\log(1/\delta)^{2\theta}+\gamma(\theta)\log(1/\delta)}{\sqrt{T}}\right)
\end{align*}
where
\[
\gamma(\theta)=\begin{cases}
64K^2 &\theta=1/2\\
8\max\left\{(4\theta)^{\theta}eK\rho,\; 4(4\theta)^{2\theta}e^2K^2\right\} & \theta\in(1/2,1]\\
\begin{split}
8\max\bigg\{
  & 2\log(2T/\delta_2)^{\theta-1}K\rho, \\ &4\left[(2^{2\theta+1}+2)\Gamma(2\theta+1)+\frac{2^{3\theta}\Gamma(3\theta+1)}{3\log(2T/\delta_2)^{\theta-1}}\right]K^2\bigg\} 
\end{split}
& \theta>1
\end{cases}
\]
and observe $\gamma(\theta) = O\left(\log(T/\delta)^{\min\{0,\theta-1\}}\right)$ for any $\theta\ge\frac12$
\end{theorem}
\begin{proof}
As with the P{\L} analysis we start the non-convex analysis with Eq.~\eqref{eq:smoothness}. From this, we get a master bound on a weighted sum of the $\|\nabla f(x_t)\|^2$. Our goal is a convergence rate for a weighted average of the $\|\nabla f(x_t)\|^2$ since this would imply convergence to a first-order stationary point, which is the best one can hope for without further assumptions. The master bound is in terms of two sums, an inner product sum and a norm sum. We bound the norm sum using an established sub-Weibull concentration inequality. The inner product sum, on the other hand, is the part that requires the MDS concentration inequality of Theorem~\ref{thm:freedman}.

\textbf{Master bound.} As in Section~\ref{sec:pl}, again assume $f$ is $L$-smooth, hence
\begin{align*}
    f(y) &\leq f(x)+\langle \nabla f(x),y-x\rangle +\frac{L}{2}\|x-y\|^2
\end{align*}
for all $x, y \in \R^d$. Set $y = x_{t+1}=x_t-\eta_t g_t=x_t-\eta_t(\nabla f(x_t)-e_t)$ and $x = x_{t}$. Then we get
\begin{align*}
    f(x_{t+1})&\leq f(x_t) -\eta_t\left(1-\frac{L\eta_t}{2}\right)\|\nabla f(x_t)\|^2+\eta_t(1-L\eta_t)\langle \nabla f(x_t),e_t\rangle+\frac{L\eta_t^2}{2}\|e_t\|^2.
\end{align*}
Summing this and using $f(x_T)\ge f\opt$, we get
\begin{align}
\label{eq:ncmaster}
    \sumt \eta_t\left(1-\frac{L\eta_t}{2}\right)\|\nabla f(x_t)\|^2 &\leq f(x_0)-f\opt+\underset{\text{inner product sum}}{\underbrace{\sumt\eta_t(1-L\eta_t)\langle \nabla f(x_t),e_t\rangle}}+\underset{\text{norm sum}}{\underbrace{\frac{L}{2}\sumt\eta_t^2\|e_t\|^2}}.
\end{align}
We would like to bound
\begin{align*}
    \sumt\eta_t(1-L\eta_t)\langle \nabla f(x_t),e_t\rangle &\leq O\left(\sumt \eta_t^2\|\nabla f(x_t)\|^2\right)\\
    \text{and }\frac{L}{2}\sumt \eta_t^2\|e_t\|^2 &\leq O\left(\sumt \eta_t^2\right)
\end{align*}
with high probability so that if $\eta_t=\Theta(1/\sqrt{t+1})$, we get
\begin{align*}
    \min_{0\leq t\leq T-1}\|\nabla f(x_t)\|^2 \leq \frac{1}{\sqrt{T}}\sumt \frac{1}{\sqrt{t+1}}\|\nabla f(x_t)\|^2 \leq O\left(\frac{\log(T+1)}{\sqrt{T}}\right)
\end{align*}
with high probability. \cmt{While these bounds are in big-$O$ notation, the bounds we prove will be precise.}

\textbf{Norm sum bound.} Assume that, conditioned on the previous iterates, $e_t$ is centered and $\|e_t\|$ is $K$-sub-Weibull($\theta$) with $\theta\ge 1/2$. Set $\eta_t=c/\sqrt{t+1}$ with $c\leq 1/L$. Let $\delta_1\in(0,1)$. Using the law of total expectation,
\begin{align*}
    \ev{\exp\left(\left(\frac{\eta_t^2\|e_t\|^2}{\eta_t^2K^2}\right)^{1/2\theta}\right)}\leq 2.
\end{align*}
Thus, $\eta_t^2$ is $\eta_t^2 K$-sub-Weibull($2\theta$) so we can apply the following sub-Weibull concentration inequality.

\begin{lemma}
\label{lma:vladimirova}
\citep[Thm. 1]{vladimirova2019sub} \citep[Lma. 5]{wong2020lasso}
Suppose $X_1,\ldots,X_n$ are sub-Weibull$(\theta)$ with respective parameters $K_1,\ldots,K_n$. Then, for all $\gamma\ge 0$,
\begin{align*}
    \pr{\bigg|\sum_{i=1}^n X_i\bigg|\ge \gamma}\leq 2\exp\left(-\left(\frac{\gamma}{v(\theta)\sum_{i=1}^n K_i}\right)^{1/\theta}\right),
\end{align*}
where $v(\theta)=(4e)^{\theta}$ for $\theta\leq 1$ and $v(\theta)=2(2e\theta)^{\theta}$ for $\theta\ge 1$.
\end{lemma}

Applying Lemma~\ref{lma:vladimirova}, we get, w.p.~$\ge 1-\delta_1$,
\begin{align*}
    \frac{L}{2}\sumt\eta_t^2\|e_t\|^2 &\leq LK^2(4e\theta\log(2/\delta_1))^{2\theta}\sumt \eta_t^2\\
    &\leq LK^2c^2(4e\theta\log(2/\delta_1))^{2\theta}\log(T+1).
\end{align*}

\textbf{Inner product sum bound.} Assume that $\theta>1$. We will prove the easier cases of $\theta=1/2$ and $\theta\in(1/2,1]$ in the appendix. Assume $f$ is $\rho$-Lipschitz continuous. Define
\begin{align*}
    \xi_t &= \eta_t(1-L\eta_t)\langle \nabla f(x_t),e_t\rangle\\
    K_{t-1} &= \eta_t(1-L\eta_t)K\|\nabla f(x_t)\|\\
    \text{and }m_t &= \eta_t(1-L\eta_t)K\rho.
\end{align*}
Recall that we defined $\F_t$ as the sigma algebra generated by $e_0,\ldots,e_t$ for all $t\ge 0$ and $\F_{-1}=\{\emptyset,\Omega\}$. So, $\xi_t$ is $\F_t$-measurable and $K_{t-1}$ is $\F_{t-1}$-measurable; hence $(\xi_t)$ and $(K_t)$ are adapted to $(\cmt{\F}_t)$. We also have, for all $t\ge 0$, $0\leq K_{t-1}\leq m_t$ almost surely, $\ev{\xi_t\mid\F_{t-1}}=0$, and
\begin{align*}
    \ev{\exp\left((|\xi_t|/K_{t-1})^{1/\theta}\right)\mid\F_{t-1}}\leq 2.
\end{align*}
In other words, $(\xi_t)$ is a sub-Weibull MDS and $(K_t)$ captures the scale parameters. Let $\delta_2,\delta_3\in(0,1)$ and define
\begin{align*}
    \delta &= \delta_2\\
    \beta &= 0\\
    \lambda &= \frac{1}{2\alpha}\\
    \text{and }x &= 2\alpha \log(1/\delta_3).
\end{align*}
Applying Theorem~\ref{thm:freedman}, we get, for all $\alpha\ge bK\rho c$, w.p.~$\ge 1-2\delta_2-\delta_3$,
\begin{align*}
    \sumt \eta_t(1-L\eta_t)\langle \nabla f(x_t),e_t\rangle &\leq 2\alpha\log(1/\delta_3)+\frac{aK^2}{\alpha}\sumt \eta_t^2(1-L\eta_t)^2\|\nabla f(x_t)\|^2.
\end{align*}
Combining this with the norm sum bound and master bound, we get, for all $\alpha\ge bK\rho c$, w.p.~$\ge 1-\delta_1-2\delta_2-\delta_3$,
\begin{align*}
    \sumt \eta_t\nu_t\|\nabla f(x_t)\|^2 \leq f(x_0)-f\opt+2\alpha\log(1/\delta_3)+ LK^2c^2(4e\theta\log(2/\delta_1))^{2\theta}\log(T+1)
\end{align*}
where
\begin{align*}
    \nu_t = 1-\frac{L\eta_t}{2}-\frac{aK^2}{\alpha}\eta_t(1-L\eta_t)^2.
\end{align*}
We want to bound $\nu_t$ away from zero. To do so, assume $c\leq \frac{1}{L}$ and $\alpha\ge 4aK^2c$. Then $\nu_t\ge \frac{1}{4}$. Setting
\begin{align*}
    \alpha = \max\{bK\rho,4aK^2\}c
\end{align*}
and plugging in $a$ and $b$ completes the proof.
\end{proof}

\begin{remark}
Note that Lipschitz continuity follows immediately if the iterates are bounded. This might lead one to consider using projected SGD, but there are certain issues preventing us from analyzing it, which we discuss in Appendices~\ref{sec:plproj} and~\ref{sec:projnc}. But, is Lipschitz continuity even a necessary assumption? \cite{lishaoji2022high}, building off of our pre-print, were actually able to relax the Lipschitz continuity assumption to the assumption that $\frac{1}{\sqrt{t+1}}\|\nabla f(x_t)\|\leq \rho~\forall t\ge 0$. To see that this works, note that we can change the definition of $m_t$ to $c(1-L\eta_t)K\rho$ and the rest of the analysis, including the result, still holds.
\end{remark}

\begin{remark}
Can we extend the analysis beyond sub-Weibull? Yes, but then we would not get logarithmic dependence on $1/\delta$. For example, if we assume $\ev{\|e_t\|^p\mid \F_{t-1}}$ for some $p>2$, then we could use Corollary 3 instead of Corollary 2 of \cite{bakhshizadeh2020sharp}. This would give us a $O(\log(1/\delta)/\sqrt{T}+1/(\delta^aT^b))$ convergence rate for some $b>1/2$. Thus, we would not have a logarithmic dependence on $1/\delta$ in general, but would approach such a dependence as the number of iterations increases.
\end{remark}

\section{Post-processing}
\label{sec:postprocessing}

Note that the results of Theorem~\ref{thm:ncconvergence} are in terms of $\frac{1}{\sqrt{T}}\sumt \frac{1}{\sqrt{t+1}}\|\nabla f(x_t)\|^2$ which is not a particularly useful quantity by itself. \cmt{To get a bound in terms of a single iterate, we prove the following probability result, introduce a novel post-processing strategy which outputs a single iterate $x$, and apply the probability result to bound $\|\nabla f(x)\|^2$ with high probability.
}

\begin{theorem}
\label{thm:validation}
Let $T\in\N$. For all $t\in[T]$, let $Z=t$ with probability $p_t$, where $\sum_{t=1}^Tp_t=1$. Let $Z_1,\ldots,Z_n$ be independent copies of $Z$. Let $Y= \{Z_1,\ldots,Z_n\}$. Let $X$ be an $\R_+^T$-valued random variable independent of $Z$. Then
\begin{align*}
    \pr{\min_{t\in Y}X_t>e\gamma} \leq \exp(-n)+\pr{\sum_{t=1}^Tp_tX_t>\gamma}~\forall \gamma>0.
\end{align*}
\end{theorem}
\begin{proof}
First, letting $\gamma>0$ and $x\in\R_+^T$,
\begin{align*}
    \pr{\min_{t\in Y}x_t>\gamma} &= \pr{\bigcap_{i=1}^n\{x_{Z_i}>\gamma\}}\\
    &\overset{(i)}{=} \prod_{i=1}^n \pr{x_{Z_i}>\gamma}\\
    &= \pr{x_Z>\gamma}^n\\
    &\overset{(ii)}{\leq} \left(\frac{1}{\gamma}\ev{x_Z}\right)^n\\
    &= \left(\frac{1}{\gamma}\sum_{t=1}^T p_tx_t\right)^n,
\end{align*}
where $(i)$ follows by the independence of the $Z_i$ and $(ii)$ follows by Markov's inequality since $x_Z$ is non-negative almost surely. Next, define
\begin{align*}
    &A = \bigg\{x\in\R_+^T~\bigg|~ \sum_{t=1}^Tp_tx_t\leq \gamma\bigg\}\\
    &B = \bigg\{(x,y)\in\R_+^T\times[T]^n~\bigg|~ x_{y_i}>e\gamma~\forall i\in[n]\bigg\}.
\end{align*}
Observe,
\begin{align*}
    \pr{(X,Y)\in B} &\overset{(i)}{=}\pr{X\in A,(X,Y)\in B}+\pr{X\in A^c,(X,Y)\in B}\\
    &\overset{(ii)}{=} \int_{A} \pr{(x,Y)\in B}\mu(dx)+\int_{A^c} \pr{(x,Y)\in B}\mu(dx)\\
    &\leq \int_{A} \left(\frac{1}{e\gamma}\sum_{t=1}^T p_tx_t\right)^n\mu(dx)+\int_{A^c} \mu(dx)\\
    &\leq \exp(-n)\int_{A}\mu(dx)+\pr{X\in A^c}\\
    &\leq \exp(-n)+\pr{\sum_{t=1}^Tp_t X_t>\gamma}
\end{align*}
where (i) follows from the law of total probability and (ii) follows from Theorem~20.3 of \cite{Billingsley} since $X$ and $Y$ are independent.
\end{proof}

\begin{corollary} \label{cor:postprocessing}
\cmt{Assume $f$ is differentiable. Let $\delta\iter\in(0,1)$ and $T\in\N$. Set $n\iter = \lceil \log(1/\delta\iter)\rceil$. Sample $n\iter$ indices with replacement from $\{0,\ldots,T-1\}$ with probabilities $p_0,\ldots,p_{T-1}$ to form the set $S=\{s_1,\ldots,s_{n\iter}\}$. Then, for $T$ iterations of SGD with any step-size sequence $(\eta_t)$,}
\begin{align*}
    \pr{\min_{t\in S}\|\nabla f(x_t)\|^2>e\gamma} \leq \pr{\sumt p_t\|\nabla f(x_t)\|^2>\gamma}+\cmt{\delta\iter}~\forall \gamma>0.
\end{align*}
\end{corollary}

\cmt{To combine the corollary with Theorem~\ref{thm:ncconvergence}, we can set}
\begin{align*}
    p_t=\frac{1/\sqrt{t+1}}{\sumt 1/\sqrt{t+1}}\text{ and use that }\frac{1}{\sumt 1/\sqrt{t+1}}\le \frac{1}{\sqrt{T}}.
\end{align*}
\cmt{To see the merit of this post-processing strategy, let's compare it to a naive approach one might take to apply the result of Theorem~\ref{thm:ncconvergence}.} The standard trick is to observe
\begin{align}
\label{eq:minsumbound}
    \min_{0\leq t\leq T-1}\|\nabla f(x_t)\|^2 \leq \frac{1}{\sqrt{T}}\sumt \frac{1}{\sqrt{t+1}}\|\nabla f(x_t)\|^2.
\end{align}
\cmt{So, we} could keep track of $\|\nabla f(x_t)\|$ at every iteration and record the iterate where \cmt{it} is lowest.  However, this requires exact gradient information, which may be more costly than the stochastic gradient used in the algorithm. In \cite{ghadimi2013stochastic}, they pick index $s$ with probability proportional to
$1/\sqrt{s+1}$ so that $\ev{\|\nabla f(x_s)\|^2}$ is proportional to the right-hand side of Eq.~\eqref{eq:minsumbound}. 
They do this for
$\Theta(\log(1/\delta))$ \textit{runs}
and pick the best of the runs. \cmt{Corollary~\ref{cor:postprocessing}}, on the other hand, \cmt{allows us to} sample a set $S$ of $n\iter=\Theta(\log(1/\delta))$ indices and pick the best \textit{iterate} from among these samples. \cmt{Hence, we call} $\delta\iter$ the iterate sampling failure probability.

\cmt{But, Corollary~\ref{cor:postprocessing} is not the end of the story since} to compute even \cmt{$\argmin_{t\in S}\|\nabla f(x_t)\|^2$} \cmt{still requires} full gradient information. In the sample average approximation setting, this can be obtained by running on the full batch of data (rather than a mini-batch). However, if this is computationally infeasible or if we are in the stochastic approximation setting, then we instead have to use empirical gradients over a test or validation set. \cmt{This is what we do for the full post-processing strategy presented in the following theorem.}

\begin{theorem} \label{thm:postprocessing}
\cmt{Let $(\Omega,\F,P)$ be a probability space. Let $F:\R^d\times \Omega\to\R$ and assume $F(\cdot;\xi)$ is differentiable for all $\xi\in\Omega$. Let $f=\ev{F(\cdot;\xi)}$.} Assume $\nabla f = \ev{\nabla F(\cdot;\xi)}$ and $\ev{\|\nabla f(x)-\nabla F(x;\xi)\|^2}\le \sigma^2~\forall x\in\R^d$. \cmt{Let $\delta\iter,\delta\emp\in(0,1)$ and $T\in\N$.
Set $n\iter = \lceil \log(1/\delta\iter)\rceil$ and} $n\emp = \lceil 6(n\iter+1)\sigma^2/(e\gamma \delta\emp)\rceil$. \cmt{Apply the following procedure:}
\begin{enumerate}
    \item \cmt{Sample $n\iter$ indices with replacement from $\{0,\ldots,T-1\}$ with probabilities $p_0,\ldots,p_{T-1}$ to form the set $S=\{s_1,\ldots,s_{n\iter}\}$.
    \item Sample $\xi_1,\ldots,\xi_{n\emp}$ independently.
    \item Run $T$ iterations of SGD with any step-size sequence $(\eta_t)$ to form $(x_t)$.}
    \item Compute \begin{align*}
    x = \underset{t\in S}{\text{argmin}}\bigg\|\frac{1}{n\emp}\sum_{j=0}^{n\emp} \nabla F(x_t;\xi_j)\bigg\|^2.
    \end{align*}
\end{enumerate}
Then,
\begin{align*}
    \pr{\|\nabla f(x)\|^2>5e\gamma}\leq \pr{\sumt p_t\|\nabla f(x_t)\|^2>\gamma}+\delta\iter+\delta\emp~\forall \gamma>0.
\end{align*}
\end{theorem}
\begin{proof}
Apply Theorem~2.4 of \cite{ghadimi2013stochastic} to get
\begin{align*}
    \pr{\|\nabla f(x)\|^2>5e\gamma}\leq\pr{\min_{t\in S}\|\nabla f(x_t)\|^2>e\gamma}+\frac{6(n\iter+1)\sigma^2}{e\gamma n\emp}~\forall \gamma,\lambda>0.
\end{align*}
\cmt{Using the definitions of $n\iter$ and $n\emp$, and applying Corollary~\ref{cor:postprocessing}, proves the result.}
\end{proof}

\cmt{We call} $\delta\emp$ the empirical gradient failure probability. \cmt{Note that if} $\|\nabla f(x)-G(x,\xi)\|$ is $K$-sub-Weibull($\theta$), \cmt{then} $e_t$ is centered and $\|e_t\|$ is $K$-sub-Weibull($\theta$), \cmt{in which case both} Theorem~\ref{thm:ncconvergence} and \cmt{Theorem~\ref{thm:postprocessing}} apply, \cmt{with} $\sigma^2=2\Gamma(2\theta+1)K^2$ (by Lemma~\ref{lma:moment}) \cmt{for the latter. Also, note that for both Corollary~\ref{cor:postprocessing} and Theorem~\ref{thm:postprocessing}, while we have to specify $T$ in advance, we only have to take $\max S$ iterations to apply the bound from Theorem~\ref{thm:ncconvergence} to the post-processing output.
}

\section{Neural Network Example}
\label{sec:numerics}

Consider the two layer neural network model
\begin{align*}
x \mapsto    \phi(x^TW)a
\end{align*}
where $\phi$ is a differentiable activation function applied coordinate-wise, $x\in\R^d$ is a data point (feature vector), $W\in\R^{d\times m}$ is the first-layer weights, and $a\in\{\pm 1\}^m$ is the second-layer weights. If we are given a data set $X\in\R^{d\times n}$ and labels $y\in\R^n$, then we can train $W$ (leaving $a$ fixed for simplicity) using the squared loss:
\begin{align*}
    f(W) = \frac{1}{2}\|\phi(X^TW)a-y\|^2.
\end{align*}
In this case,
\begin{align*}
    \text{vec}(\nabla f(W)) = \left(\diag(a)\phi'(W^TX)* X\right)\left(\phi(X^TW)a-y\right)
\end{align*}
where $*$ is the Khatri-Rao product \citep{oymak2020toward}. \cmt{For} our example, we use \cmt{the} GELU activation, $\phi(x)=x[1+\text{erf}(x/\sqrt{2})]/2$ \citep{hendrycks2016gaussian}, \cmt{which satisfies $|\phi'(x)|\le 1.13$ and $|\phi''(x)|\le .11$ for all $x\in\R$. Thus, we can apply Lemma~\ref{lma:oymak}. Since $\lim_{x\to\infty}\phi''(x)=0=\lim_{x\to-\infty}\phi''(x)$, we heuristically set $b=0$ in the lemma and estimate the strong smoothness of $f$ as} $\approx m\|X\|_2^2$, \cmt{setting our step-size accordingly.}

\begin{figure}
    \centering
    \includegraphics[width=.95\linewidth]{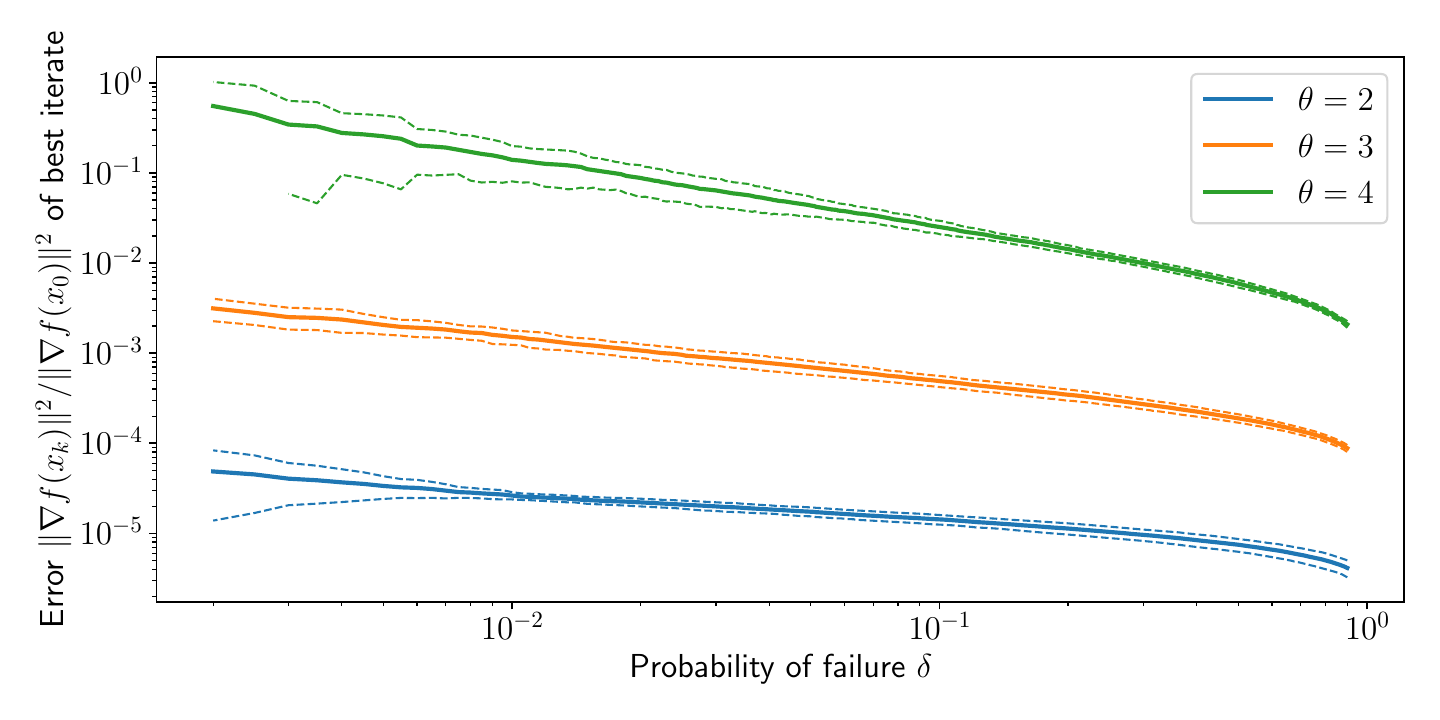}
    \caption{Empirical $1-\delta$ convergence error, averaged over 10000 runs. 
    The dashed lines show the mean $\pm$ one standard deviation (computed over 5 blocks of 2000 runs each).
    The data are less reliable for small $\delta$.
    }
    \label{fig:together_ver2}
\end{figure}

In order to demonstrate the effect of gradient noise on convergence error, we train $W$ on a fixed synthetic data set while injecting noise into the gradient. The labels come from a neural network model with width, $m'$, larger than the width of the training model. We also make sure the total number of trainable parameters is less than $n$ so that $f\opt>0$. The noise we inject has uniformly random direction and Weibull norm with scale parameter $K=1$ and shape parameter $1/\theta$, with $\theta\in\{2,3,4\}$. 
We keep the same initialization for all trials. We run $100$ iterations of SGD and then define the convergence error to be the best gradient norm squared divided by the initial gradient norm squared. We compute the empirical CDF of the convergence error over 10000 trials and then consider $\delta$ in the ranges $[0.1,0.2]$, $[0.01,0.1]$, and $[0.001,0.01]$. We care about the dependence of the convergence error on $\delta$ for small $\delta$, but for too small of $\delta$, the empirical CDF is not a good approximation to the true CDF (see Fig.~\ref{fig:together_ver2}) due to the nature of order statistics. Our code can be found at \url{https://github.com/liammadden/sgd}.

From Theorem~\ref{thm:ncconvergence}, for a particular range of $\delta$ (and for fixed $T$), the upper bound has dependence either $\log(1/\delta)^{\theta}$, $\log(1/\delta)$, or $\log(1/\delta)^{2\theta}$. For sufficiently small $\delta$, the $\log(1/\delta)^{2\theta}$ will dominate, but the upper limit of this range of $\delta$ may be smaller than the lower limit of the range of $\delta$ for which the empirical CDF is a good approximation to the true CDF. In other words, if we showed that the true CDF for this particular example has $\log(1/\delta)^{2\theta}$ dependence for $\delta$ sufficiently small, then we would have shown that the $\delta$-dependence of the upper bound in Theorem~\ref{thm:ncconvergence} is tight. However, with the empirical CDF, we are only able to show that the exponent increases as $\delta$ shrinks. 
Figure~\ref{fig:together_ver2} shows
the dependence for the three different ranges. We assume the convergence error has dependence $b\log(1/\delta)^a$ and find the line of best fit as $\log(b)+a\log\log(1/\delta)$. 
The different values of $a$ are given in Table~\ref{tab:slopes}. 

Our experiments suggest that injected Weibull noise results in convergence error with dependence $\Omega(\log(1/\delta)^{c(\delta)\theta+d(\delta)})$ where $c(\delta)$ increases as $\delta$ decreases, thus roughly corroborating our upper bound. In particular, by computing a line of best fit for the $a$ values in each of the three $\delta$ ranges, we can estimate that $c(\delta)$ increases from 0.43, to 0.77, to 1.34 and $d(\delta)$ decreases from -0.28, to -0.7, to -1.42, suggesting that we are in the $\log(1/\delta)^{\theta-1}$ regime.

\begin{table}
    \centering
    \begin{tabular}{cccc}
    \toprule
    $a$ & $\delta\in[0.1,0.2]$ & $\delta\in[0.01,0.1]$ & $\delta\in[0.001,0.01]$\\
    \midrule
    $\theta=2$ & 0.64 & 0.85 & 1.40\\
    \midrule
    $\theta=3$ & 0.93 & 1.56 & 2.33\\
    \midrule
    $\theta=4$ & 1.51 & 2.38 & 4.08 \\
    \bottomrule
    \end{tabular}
    \caption{Empirically estimated exponents of $\log(1/\delta)$}
    \label{tab:slopes}
\end{table}

\begin{figure}
    \centering
    \includegraphics[width=.95\linewidth]{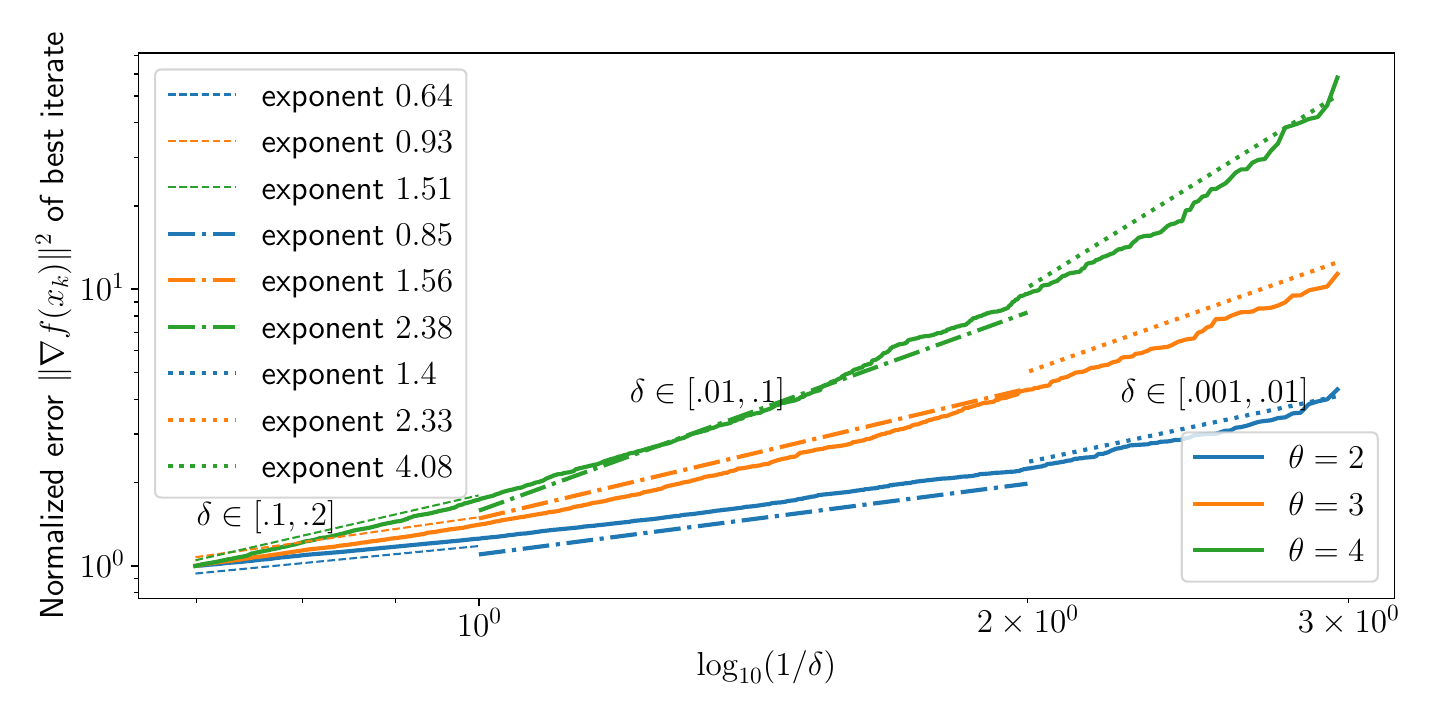}
    \caption{Same data as Fig.~\ref{fig:together_ver2} but each line series is normalized, and the x-axis is $\log(1/\delta)$ \emph{and} plotted on a logarithmic scale, so $\log(1/\delta)^a$ dependence shows us a straight line with slope $a$. The $\delta$ range is from $0.2$ (left side) to $0.01$ (right side), since any smaller $\delta$ has unreliable statistics.
    Lines of best fit using the exponents from Table~\ref{tab:slopes} are shown (with arbitrary shifts for clarity).
    }
    \label{fig:together_with_fits}
\end{figure}

\section{Conclusions}

This paper analyzed the convergence of SGD for objective functions that satisfy the P{\L} condition and for generic non-convex objectives. Under a sub-Gaussian assumption on the gradient error, we showed 
a high probability convergence rate matching the mean convergence rate for P{\L} objectives. Under a sub-Weibull assumption on the noise, we showed a high probability convergence rate matching the mean convergence rate for non-convex objectives. We also provided and analyzed a post-processing method for choosing a single iterate. To prove the convergence rate, we first proved a Freedman-type inequality for martingale difference sequences that extends previous Freedman-type inequalities beyond the sub-exponential threshold to allow for sub-Weibull tail-decay. Finally, we considered a synthetic neural net problem and showed that the heaviness of the tail of the gradient noise has a direct effect on the heaviness of the tail of the convergence error.


\acks{All three authors gratefully acknowledge support from the National Science Foundation (NSF), Division of Mathematical Sciences, through the award \# 1923298.}


\vskip 0.2in

\appendix

\cmt{
\section{Proof of Lemma~\ref{lma:assump}}
\label{sec:assumplmapf}

Here we prove Lemma~\ref{lma:assump}.

\begin{proof}
First,
\begin{align*}
    f(W) = \sum_{i=1}^n \frac{1}{2}\left[\phi(x_i^\top W)v-y_i\right]^2
\end{align*}
and
\begin{align*}
    \frac{\partial}{\partial w_{s,t}}\left[\phi(x_i^\top W)v-y_i\right] &= \frac{\partial}{\partial w_{s,t}}\sum_{j=1}^m v_j\phi(x_i^\top w_j)= v_t\phi'(x_i^Tw_t)x_{s,i}.
\end{align*}
So,
\begin{align*}
    \frac{\partial}{\partial w_{s,t}}f(W) = \sum_{i=1}^n\left[\phi(x_i^\top W)v-y_i\right]v_t\phi'(x_i^Tw_t)x_{s,i}.
\end{align*}
Thus,
\begin{align*}
    \frac{\partial^2}{\partial w_{\ell,r}\partial w_{s,t}}f(W) &= \sum_{i=1}^nv_r\phi'(x_i^\top w_r)x_{\ell,i}v_t\phi'(x_i^\top w_t)x_{s,i}\\
    &\quad+\delta_{r,t}\sum_{i=1}^n\left[\phi(x_i^\top W)v-y_i\right]v_t\phi''(x_i^\top w_t)x_{\ell,i}x_{s,i}
\end{align*}
where $\delta$ is the Kronecker delta. Let $b$ be the largest absolute element of $X$. Then, viewing $W$ as a vector,
\begin{align*}
    \|\nabla f(W)\|_2 \le \sqrt{md}\|\nabla f(W)\|_{\infty}\le \sqrt{md}nab\|v\|_{\infty}(a\sqrt{m}\|v\|_2+\|y\|_{\infty})
\end{align*}
and
\begin{align*}
    \|\nabla^2 f(W)\|_2 \le md\|\nabla^2 f(W)\|_{\max}\le mdn\left(a^2b^2\|v\|_{\infty}^2+ab^2\|v\|_{\infty}(a\sqrt{m}\|v\|_2+\|y\|_{\infty})\right),
\end{align*}
proving the result.
\end{proof}

And here we prove Lemma~\ref{lma:oymak}.

\begin{proof}
Define $F:w\mapsto \phi(X^\top\text{vec}^{-1}(w))v$ and $\ls:w\mapsto \|F(w)-y\|_2^2/2$. Then $\|DF(w)\|_2\le \sqrt{m}a\|v\|_{\infty}\|X\|_2$ and $\|DF(w)-DF(u)\|_2\le b\|v\|_{\infty}\|X\|_2\|X\|_{1,2}\|w-u\|_2$ for all $w,u\in\R^{md}$ by Lemmas 3 and 5 of \cite{oymak2020toward}. Let $\rho=\sqrt{m}a\|v\|_{\infty}\|X\|_2$ and $L=b\|v\|_{\infty}\|X\|_2\|X\|_{1,2}$. Observe
\begin{align*}
    \|\nabla \ls(w)\|_2 &= \|DF(w)^\top(F(w)-y)\|_2\le \|DF(w)\|_2\|F(w)-y\|_2\le \rho \sqrt{2\alpha}
\end{align*}
and
\begin{align*}
    \|\nabla \ls(w)-\nabla \ls(u)\| &= \|DF(w)^\top(F(w)-F(v))\|_2+\|(DF(w)-DF(v))^\top(F(v)-y)\|_2\\
    &\le \|DF(w)\|_2\|F(w)-F(v)\|_2+\|DF(w)-DF(v)\|_2\|F(v)-y\|_2\\
    &\le \rho^2\|w-v\|_2+L\|w-v\|_2\sqrt{2\alpha}=(\rho^2+L\sqrt{2\alpha})\|w-v\|_2,
\end{align*}
proving the result.
\end{proof}
}

\section{Remaining Work for Proof of Theorem~\ref{thm:freedman}}
\label{sec:freedmanpf}

First we need the following two lemmas. Lemma~\ref{lma:subexp} extends Proposition~2.7.1 of \cite{Vershynin} to interpolate between the sub-Gaussian and sub-exponential regimes.

\begin{lemma}
\label{lma:subgau}
\citep[Prop. 2.5.2(e)]{Vershynin}
If $X$ is centered and $K$-sub-Guassian then $\ev{\exp(\lambda X)}\leq \exp(\lambda^2 K^2)~\forall \lambda\in\R$.
\end{lemma}

\begin{lemma}
\label{lma:subexp}
If $X$ is centered and $K$-sub-Weibull($\theta$) with $\theta\in(1/2,1]$, then\\ $\ev{\exp(\lambda X)}\leq\exp\left(\frac{\lambda^2}{2}(4\theta)^{2\theta}e^2K^2\right)$ for all $\lambda\in\left[0,\frac{1}{(4\theta)^{\theta}eK}\right]$.
\end{lemma}
\begin{proof}
First, using Lemma~\ref{lma:moment} and $\Gamma(x+1)\leq x^x~\forall x\ge 1$, we can get $\|X\|_p\leq (2\theta)^{\theta}Kp^{\theta}$ for all $p\ge 1/\theta$, and so, in particular, for all $p\ge 2$. Thus, if $\lambda\in\left[0,\frac{1}{(4\theta)^{\theta}eK}\right]$, then
\begin{align*}
    \ev{\exp(\lambda X)} &= \ev{1+\lambda X+\sum_{p=2}^{\infty}\frac{(\lambda X)^p}{p!}}\\
    &= 1+\sum_{p=2}^{\infty}\frac{\lambda^p\|X\|_p^p}{p!}\\
    &\leq 1+\sum_{p=2}^{\infty}\frac{\lambda^p(2\theta)^{\theta p}K^pp^{\theta p}}{p!}\\
    &\leq 1+\sum_{p=2}^{\infty}\left(\frac{\lambda(2\theta)^{\theta}eK}{p^{1-\theta}}\right)^p\\
    &\leq 1+\sum_{p=2}^{\infty}\left(\lambda(4\theta)^{\theta}(e/2)K\right)^p\\
    &\leq 1+\frac{\left(\lambda(4\theta)^{\theta}(e/2)K\right)^2}{1-\lambda(4\theta)^{\theta}(e/2)K}\\
    &\leq 1+2\left(\lambda(4\theta)^{\theta}(e/2)K\right)^2\\
    &\leq \exp\left(\frac{\lambda^2}{2}(4\theta)^{2\theta}e^2K^2\right),
\end{align*}
completing the proof.
\end{proof}

Then, the next three lemmas allow us to include previous results as special cases of the theorem.

\begin{lemma}
\label{lma:freedmanfan}
\cite[Thm.~2.6]{Fan}
Let $(\Omega,\F,(\F_i),\P)$ be a filtered probability space.
Let $(\xi_i)$ and $(K_i)$ be adapted to $(\F_i)$. Let $n\in\N$. For all $i\in[n]$, assume $0\leq K_{i-1}\leq m_i$ almost surely, $\ev{\xi_i\mid\F_{i-1}}=0$, and
\begin{align*}
    \ev{\exp(\lambda \xi_i)\mid\F_{i-1}}\leq \exp\left(\frac{\lambda^2}{2}aK_{i-1}^2\right)~\forall\lambda\in\left[0,\frac{1}{bK_{i-1}}\right].
\end{align*}
Then, for all $x,\beta\ge 0$, and $\lambda\in\left[0,\frac{1}{b\max_{i\in[r]}m_i}\right]$,
\begin{align*}
    \pr{\bigcup_{k\in[n]}\bigg\{\sum_{i=1}^k\xi_i\ge x\text{ and }\sum_{i=1}^k aK_{i-1}^2\leq \beta\bigg\}}&\leq \exp\left(-\lambda x+\frac{\lambda^2}{2}\beta \right).
\end{align*}
\end{lemma}
\begin{proof}
Define
\begin{align*}
    \psi_i = \exp\left(\lambda \xi_i\right)
\end{align*}
and
\begin{align*}
    A_k =\bigg\{\sum_{i=1}^k \xi_i\ge x\text{ and }\sum_{i=1}^k aK_{i-1}^2\leq \beta\bigg\}.
\end{align*}
Then $\omega\in A_k$ implies
\begin{align*}
    \prod_{i=1}^k \frac{\ev{\psi_i\mid\F_{i-1}}}{\psi_i} &\leq \exp\left(-\lambda\sum_{i=1}^k \xi_i+\frac{\lambda^2}{2}\sum_{i=1}^k aK_{i-1}^2\right)\\
    &\leq \exp\left(-\lambda x+\frac{\lambda^2}{2}\beta \right).
\end{align*}
\end{proof}

\begin{lemma}
\label{lma:freedmanharvey}
\cite[Thm.~3.3]{Harvey1}
Let $(\Omega,\F,(\F_i),\P)$ be a filtered probability space.
Let $(\xi_i)$ and $(K_i)$ be adapted to $(\F_i)$. Let $n\in\N$. For all $i\in[n]$, assume $K_{i-1}\ge 0$ almost surely, $\ev{\xi_i\mid\F_{i-1}}=0$, and
\begin{align*}
    \ev{\exp(\lambda \xi_i)\mid\F_{i-1}}\leq \exp\left(\frac{\lambda^2}{2}aK_{i-1}^2\right)~\forall\lambda\ge 0.
\end{align*}
Then, for all $x,\beta,\alpha \ge 0$, and $\lambda\in\left[0,\frac{1}{2\alpha}\right]$,
\begin{align*}
    \pr{\bigcup_{k\in[n]}\bigg\{\sum_{i=1}^k\xi_i\ge x\text{ and }\sum_{i=1}^k aK_{i-1}^2\leq \alpha\sum_{i=1}^k \xi_i+\beta\bigg\}}&\leq \exp(-\lambda x+2 \lambda^2\beta).
\end{align*}
\end{lemma}

\begin{lemma}
\label{lma:freedman}
\citep{Freedman}
Let $(\Omega,\F,(\F_i),\P)$ be a filtered probability space.
Let $(\xi_i)$ and $(K_i)$ be adapted to $(\F_i)$. For all $i\in[n]$, assume $K_{i-1}\ge 0$ almost surely, $\ev{\xi_i\mid\F_{i-1}}=0$, and
\begin{align*}
    \ev{\exp(\lambda \xi_i)\mid\F_{i-1}}\leq \exp\left(\frac{\lambda^2}{2}aK_{i-1}^2\right)~\forall\lambda\ge 0.
\end{align*}
Then, for all $x,\beta\ge 0$, and $\lambda\ge 0$,
\begin{align*}
    \pr{\bigcup_{k\in[n]}\bigg\{\sum_{i=1}^k\xi_i\ge x\text{ and }\sum_{i=1}^k aK_{i-1}^2\leq \beta\bigg\}}&\leq \exp\left(-\lambda x+\frac{\lambda^2}{2}\beta \right).
\end{align*}
\end{lemma}

If the $\xi_i$'s are sub-Gaussian, that is, if $\theta=1/2$, then, from Lemma~\ref{lma:subgau},
\begin{align*}
    \ev{\exp(\lambda\xi_i)\mid\F_{i-1}}\leq\exp\left(\frac{\lambda^2}{2}2K_{i-1}^2\right)~\forall\lambda\in\R,
\end{align*}
so we can apply Lemma~\ref{lma:freedmanharvey} if $\alpha>0$ or Lemma~\ref{lma:freedman} if $\alpha=0$.

If the $\xi_i$'s are at most sub-exponential, that is, if $1/2<\theta\leq 1$, then, from Lemma~\ref{lma:subexp},
\begin{align*}
    \ev{\exp(\lambda\xi_i)\mid\F_{i-1}}\leq\exp\left(\frac{\lambda^2}{2}(4\theta)^{2\theta}e^2K_{i-1}^2\right)~\forall\lambda\in\left[0,\frac{1}{(4\theta)^{\theta}eK_{i-1}}\right],
\end{align*}
so we can apply Lemma~\ref{lma:subexpharvey} if $\alpha>0$ or Lemma~\ref{lma:freedmanfan} if $\alpha=0$.

\section{P{\L} Projected SGD}
\label{sec:plproj}

When optimizing over a constraint set, $X \subsetneq \R^d$, if $f$ is strongly convex, then so is $f$ plus the indicator function of $X$, and results for gradient descent methods easily extend to projected gradient descent. On the other hand, if $f$ is P{\L}, then $f$ plus the indicator function is not KL (Kurdyka-Łojasiewicz). This has real impacts on gradient descent algorithms, where gradient descent might converge while projected gradient descent does not. For example, there is a smooth function, a mollified version of $f(x,y)=\left(a(x)_+^2-b\left(|y|-c\right)_+\right)_+$, such that the P{\L} inequality is satisfied but projected gradient descent does \textit{not} converge to a minimizer;  we formalize this in the remark below.

\begin{figure}
    \centering
    \includegraphics[width=.6\linewidth]{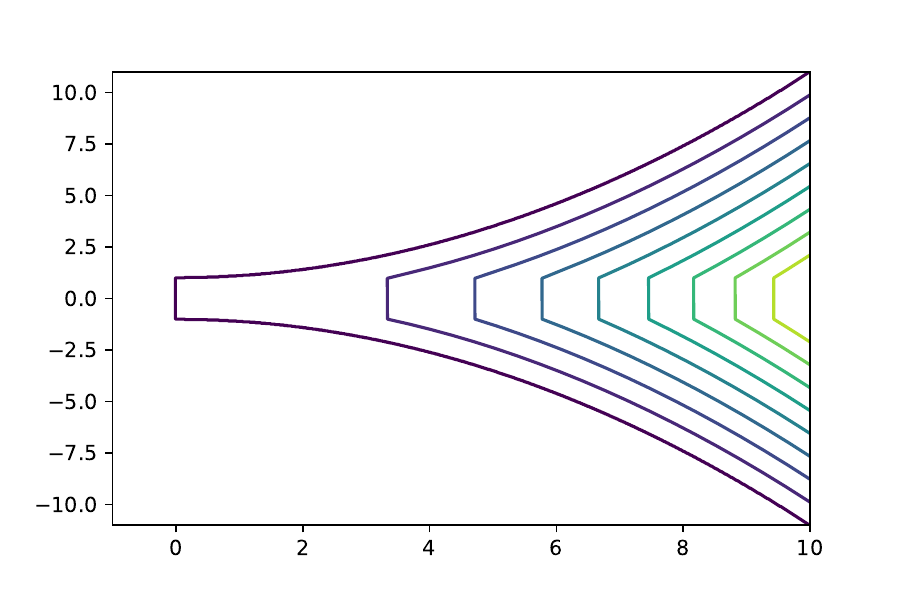}
    \caption{Contour plot of the P{\L} function counter-example to projected gradient flow}
    \label{fig:projectedpl}
\end{figure}

\begin{remark}
\label{rmk:pl}
Consider $f(x,y)=\left(a(x)_+^2-b\left(|y|-c\right)_+\right)_+$ where $(\cdot)_+$ denotes $\max(\cdot,0)$ and $a,b>0,c\geq 0$. The minimum of $f$ is 0 and $X\opt=\{(x,y)~|~x\leq 0 \text{ or } |y|\geq \frac{a}{b}x^2+c\}$. If we use $\varphi(x)=\X_{B_1(0)}\cdot\exp\left(-1/(1-\|x\|^2)\right)/\Phi$---where $\X$ denotes the indicator function, $B_1(0)$ denotes the ball of radius 1 centered at 0, and $\Phi$ is the normalization constant---to mollify $f$, then, for $\epsilon<c$, $f_{\epsilon}$ has P{\L} constant $2a$ and smoothness constant $2a$. Consider the starting point $(d,0)$. For $a,b,c,d$ chosen appropriately, the distance from $(d,0)$ to its projection onto $X\opt$ is strictly less than $d$. Thus, if we let $X$ be the ball centered at $(d,0)$ with radius equal to exactly that distance, then the constrained problem and the unconstrained problem have the same minimum. However, projected gradient flow, starting from $(d,0)$, ends up stuck at a non-minimizer: the point of $X$ closest to $(0,0)$. See Figure~\ref{fig:projectedpl} for the contour plot when $a=1/10$, $b=1=c$, and $d=10$.
\end{remark}

In order to generalize gradient methods to projected gradient methods under PL-like assumptions, the proper generalization is that the function should satisfy a \emph{proximal} P{\L} inequality~\citep{Karimi}. However, such an assumption is quite restrictive compared to the P{\L} inequality. We leave the problem of convergence with just the P{\L} inequality, via added noise or a Frank-Wolfe type construction, as a future research direction.

\section{Non-convex Projected SGD}
\label{sec:projnc}

Consider $x_{t+1}=\proj(x_t-\eta_tg_t)=\proj(x_t-\eta_t(\nabla f(x_t)-e_t))$. Define
\begin{align*}
    &G_t = \frac{x_t-\proj(x_t-\eta_t\nabla f(x_t))}{\eta_t}\\
    &E_t = \frac{x_{t+1}-\proj(x_t-\eta_t\nabla f(x_t))}{\eta_t}.
\end{align*}
Note that if $\proj=I$, then $G_t=\nabla f(x_t)$ and $E_t=e_t$. Moreover, $x_t=\proj (x_t)$ and $x_{t+1}=\proj (x_{t+1})$ so, by the non-expansiveness of $\proj$, $\|G_t\|\leq \|\nabla f(x_t)\|$ and $\|E_t\|\leq \|e_t\|$. We can get even tighter bounds using the second prox theorem \cite[Thm. 6.39]{Beck}: $\|G_t\|^2\leq \langle \nabla f(x_t),G_t\rangle$ and $\langle G_t,E_t\rangle \leq \langle \nabla f(x_t),E_t\rangle$.

It is easy to come up with an example where $\|\nabla f(x_t)\|$ does not go to zero, so we would like to bound $\|G_t\|$ instead. We start as usual:
\begin{align*}
    f(x_{t+1}) \leq f(x_t)+\langle \nabla f(x_t),x_{t+1}-x_t\rangle +\frac{L}{2}\|x_{t+1}-x_t\|^2.
\end{align*}
Focusing on the norm term,
\begin{align*}
    \frac{L}{2}\|x_{t+1}-x_t\|^2 &= \frac{L\eta_t^2}{2}\|G_t\|^2-L\eta_t^2\langle G_t,E_t\rangle +\frac{L\eta_t^2}{2}\|E_t\|^2.
\end{align*}
Focusing on the inner product term,
\begin{align*}
    \langle \nabla f(x_t),x_{t+1}-x_t\rangle &= \eta_t\langle \nabla f(x_t),E_t-G_t\rangle\\
    &= \eta_t\langle \nabla f(x_t),E_t\rangle-\eta_t\langle \nabla f(x_t),G_t\rangle\\
    &\nleq \eta_t\langle G_t,E_t\rangle-\eta_t\|G_t\|^2.
\end{align*}
Unfortunately, we cannot proceed any further. \cite{ghadimi2016mini} are able to get around this but at the cost of getting $\sumt \eta_t\|e_t\|^2=O(\sqrt{T})$ instead of $\sumt \eta_t^2\|e_t\|^2=O(\log(T))$. To mitigate this, they require an \textit{increasing} batch-size. \cite{reddi2016proximal} were able to remove this requirement, but only for non-convex projected SVRG \textit{not} non-convex projected SGD. Thus, we leave the analysis of non-convex projected SGD as an open problem.

\vskip 0.2in
\bibliography{main}

\end{document}